\documentclass[12pt,a4paper]{amsart}
\usepackage[utf8]{inputenc}	
\usepackage[T1]{fontenc}
\usepackage[in,headings]{fullpage}
\usepackage{amsbsy, amscd, amsfonts, amsmath, amsrefs, amssymb, amsthm}
\usepackage{graphicx}
\usepackage{float}
\usepackage{color}   
\usepackage[colorlinks=true,linkcolor=blue,citecolor=blue,urlcolor=blue]{hyperref}

\newtheorem{theorem}{Theorem}%[section]

\newtheorem{corollary}[theorem]{Corollary}
\newtheorem{proposition}[theorem]{Proposition}
\newtheorem{lemma}[theorem]{Lemma}

\theoremstyle{definition}
\newtheorem{definition}[theorem]{Definition}
\newtheorem{remark}[theorem]{Remark}

\theoremstyle{remark}

\newcommand{\C}{\mathbf{C}}

\newcommand{\R}{\mathbf{R}}

\renewcommand{\Re}{\mathop{\mathrm{Re}}\nolimits}
\renewcommand{\Im}{\mathop{\mathrm{Im}}\nolimits}
\newcommand{\Rzeta}{\mathop{\mathcal R }\nolimits}

\newfont{\cmbsy}{cmbsy10}
\newfont{\cmmib}{cmmib10}
\newcommand{\Orden}{\mathop{\hbox{\cmbsy O}}\nolimits}

\begin{document}
\title{Asymptotic Expansions of the auxiliary function}
\author[Arias de Reyna]{J. Arias de Reyna}
\address{%
Universidad de Sevilla \\ 
Facultad de Matem\'aticas \\ 
c/Tarfia, sn \\ 
41012-Sevilla \\ 
Spain.} 

% AMS subject classifications (used in AMS journals)

% AMS keywords (used in AMS journals)

% acknowledge support, etc
%\thanks{This research was  supported by MINECO grant MTM2015--63699-P}
% \thanks{We would like to thank our colleagues for their helpful
%  criticism.}

\email{arias@us.es, ariasdereyna1947@gmail.com}

%\date{December 12, 2023, \texttt{100-Asymptotic-v7.tex}}

\begin{abstract}
Siegel in 1932 published a paper on Riemann's posthumous writings, including a study of the Riemann-Siegel formula. 
In this paper we explicitly give the asymptotic developments of $\Rzeta(s)$ suggested by Siegel. We extend the range of validity of these asymptotic developments.  As a consequence we specify a region in which the function $\Rzeta(s)$ has no zeros. We also give complete proofs of some of Siegel's assertions.

We also include a theorem on the asymptotic behaviour of $\Rzeta(\frac12-it)$ for $t \to+\infty$. Although the real part of $e^{-i\vartheta(t)}\Rzeta(\frac12-it)$ is $Z(t)$ the imaginary part grows exponentially, this is why for the study of the zeros of $Z(t)$ it is preferable to consider $\Rzeta(\frac12+it)$ for $t>0$.
\end{abstract}

\maketitle

\section{Introduction.}
We study the asymptotic expansions of $\Rzeta(s)$ considered in \cite{Siegel}. We improve on Siegel in several points. Extending the domain of validity  of  
asymptotic expansions.  Giving the form of the terms and obtaining them in 
a way that, at most, is only suggested by Siegel.  We also give a complete proof of equation  (85) of Siegel, extending and specifying his statements. 

We consider here the original Riemann-Siegel expansion in \cite{Siegel}. This is different from the one given in Arias de Reyna \cite{A86}. The terms of the Riemann-Siegel expansion in \cite{Siegel} are analytic functions. The terms in the expansion in \cite{A86} are not analytic functions. This has several advantages that I will consider in a forthcoming paper. The main difference is that $a=\sqrt{t/2\pi}$ in \cite{A86} that is always real, is replaced by $\eta=\sqrt{(s-1)/2\pi i}$, that is complex. But the coefficients are related; for example, our polynomials $P_k(z)$ in \eqref{gexpansion} coincide with the polynomials $P_k(z)$ in \cite{A86}*{eq.~(3.7)} when in \cite{A86} we take $\sigma=0$. 

The derivation of the expansion in \cite{Siegel} starts from the first representation integral of $\zeta(s)$ in the Riemann paper. While Gabcke \cite{G}, and I following him \cite{A86}, start from the integral representation of $\Rzeta(s)$. Here I also follow this path, so that our derivation of Siegel expansions in \cite{Siegel} would be simpler and very similar to \cite{A86}.

The main asymptotic expansion for $\Rzeta(s)$ is given in Theorem \ref{Theorem1}, it is valid on the plane except for any small angle with vertex at $0$ and bisector the negative real axis. 
Using the relation $\zeta(s)=\Rzeta(s)+\chi(s)\overline{\Rzeta}(1-s)$ another asymptotic expansion is found in Theorem \ref{P:intermediate} valid on the plane except an angle with vertex at $1$ and bisector the half line $[1,+\infty)$.

From \ref{Theorem1} we derive in Theorem \ref{theorem6} that $\Rzeta(s)$ is approximated by an adequate zeta sum in the first and part of the fourth quadrant. In Theorem \ref{asympinf} obtains an asymptotic expression valid in most of the third and part of the fourth quadrant. There is a narrow zone in the fourth quadrant separating the regions in which the above Theorems apply. In this zone  a line of zeros of the function $\Rzeta(s)$ is found. 

From Theorem \ref{P:intermediate} we derive Theorem \ref{T: thm N}, an asymptotic expression valid for the third quadrant and part of the second.  This proves that in this region the zeros are the trivial ones at $s=-2n$ for $n=1$, $2$, \dots

The behaviour of $\Rzeta(s)$ for $t>0$ and $-2<-\sigma \le  t^\alpha$ is delicate. Siegel proved that it contains many zeros of $\Rzeta(s)$ connected with the zeros of $\zeta(s)$.
Following Siegel, we  consider the behaviour in this zone in Theorem  \ref{left}, and Corollaries \ref{oldcor}, \ref{newcor}  and \ref{cor85-2}.

For $s=\frac12+it$, the Riemann-Siegel function is given by $Z(t)=2\Re\{\pi^{-s/2}\Gamma(s/2)\Rzeta(s)\}$. This is true for real $t$, Siegel only considers the case $t>0$. For good reason. We  give (see Theorem \ref{asymptinf}) the asymptotic value of $\Rzeta(\tfrac12 -it)$ for $t\to+\infty$, which shows why there is not much interest in this equation for $t<0$.

\section{Formal derivation of the first expansion.}\label{secdos}

Start from the definition 

\begin{equation}\label{dos}
\Rzeta(s)=\int_{0\swarrow1}\frac{x^{-s}
e^{\pi i x^2}}{e^{\pi i x}-e^{-\pi i x}}\,dx.
\end{equation}
The integrand is meromorphic for $x$ in the complex plane with a cut along the negative real axis, taking $x^{-s}=\exp(-s\log x)$ with $\log x$ denoting the principal value. 
The saddle point for this integral is at any value of the square root $\xi=\sqrt{s/2\pi i}$. There are two possible values, but we take the root contained in the sector
$-3\pi/4<\arg\xi<\pi/4$ so that we can change the path of integration to a parallel line passing through $\xi$ within the domain of the integrand. 

In what follows, we assume that $s$ is given with $-\pi<\arg(s)<\pi$. So $s\in\C$
is not a negative real number, but except for this restriction may be arbitrary.
Then we define
\begin{gather}\label{defeta}
\xi:=\sqrt{\frac{s}{2\pi i}},\quad -\frac{3\pi}{4}<\arg\xi<\frac{\pi}{4},
\quad\xi_1:=\Re(\xi),\quad \xi_2:=\Im(\xi),\\
\ell:=\lfloor\xi_1-\xi_2\rfloor,\quad w(z):=\exp\Bigl\{-2\pi i\xi^2\Bigl(
\log\Bigl(1+\frac{z}{\xi}\Bigr)-\frac{z}{\xi}+\frac12\Bigl(\frac{z}{\xi}
\Bigr)^2\Bigr)\Bigr\}.\label{defm}
\end{gather}
We transform the integrand in the usual way (see \cite{G} or \cite{A86}). 
The integrand is a meromorphic function in the region $-\pi<\arg x<\pi$.
Due to the election of the argument of $\xi$, we have $\ell\ge0$. So we may apply
Cauchy's Theorem to move the line of integration between $\ell$ and $\ell+1$ 
without leaving the region $-\pi<\arg x<\pi$.
\begin{equation}\label{Cauchy}
\Rzeta(s)=\sum_{n=1}^\ell n^{-s}+\xi^{-s}e^{\pi i \xi^2}
\int_{\ell\swarrow \ell+1}\frac{e^{2\pi i (x-\xi)^2}}
{e^{\pi i x}-e^{-\pi i x}}w(x-\xi)\,dx
\end{equation}
Now we change variable (from $x$ to $v$) in the integral, where
\begin{equation}\label{defp}
x=\frac{v}{2}+\ell+\frac12,\quad q:=-2(\ell+\tfrac12 -\xi),\quad x-\xi=\frac12(v-q)
\end{equation}
and introduce
\begin{equation}
z:=i\sqrt{\pi}(v-q),\qquad \tau:=-\frac{1}{4\sqrt{\pi}\,\xi}.
\end{equation}
Then we have
\begin{multline*}
2\pi i \xi^2\Bigl(\log\Bigl(1+\frac{x-\xi}{\xi}\Bigr)-\frac{x-\xi}{\xi}
+\frac12\Bigl(\frac{x-\xi}{\xi}\Bigr)^2\Bigr)\\
=2\pi i \xi^2\log\Bigl(1+\frac{v-q}{2\xi}\Bigr)-
2\pi i \xi \frac{v-q}{2}+\pi i \xi^2\frac{(v-q)^2}{4\xi^2}
=\frac{i}{8\tau^2}\log(1+2i\tau z)+\frac{z}{4\tau}-i\frac{z^2}{4},
\end{multline*}
$2\pi i(x-\xi)^2=\frac{\pi i}{2}(v-q)^2$ and
\[e^{\pi i x}-e^{-\pi i x}=2i\sin\pi x=2i\sin\pi\Bigl(\frac{v}{2}+\ell+\frac12\Bigr)
=2i\cos\pi\Bigl(\frac{v}{2}+\ell\Bigr)=2i(-1)^\ell\cos\frac{\pi v}{2}.\]
Equation \eqref{Cauchy} transform into
\begin{equation}\label{Cauchy2}
\Rzeta(s)=\sum_{n=1}^\ell n^{-s}+\frac{(-1)^\ell}{2i}\xi^{-s}
e^{\pi i \xi^2}\int_{0\swarrow}\frac{e^{\frac{\pi i}{2}(v-q)^2}}
{2\cos\frac{\pi v}{2}}g\bigl(\tau,i\sqrt{\pi}(v-q)\bigr)\,dv
\end{equation}
where the transformed path of integration pass through $0$ and also near $q$ and
\begin{equation}\label{defg}
g(\tau,z):=\exp\Bigl\{-\frac{i}{8\tau^2}\log(1+2i\tau z)-\frac{z}{4\tau}
+i\frac{z^2}{4}\Bigr\}.
\end{equation}

For any non negative integer $K$ we consider the expansion
\begin{equation}\label{gexpansion}
g(\tau,z)=\sum_{k=0}^\infty P_k(z)\tau^k=\sum_{k=0}^KP_k(z)\tau^k+Rg_K(\tau,z),
\qquad (2|\tau z|<1).
\end{equation}
This defines $Rg_K(\tau,z)$  whenever $g(\tau,z)$ is defined, independently of the convergence of the Taylor series.  
We insert this expansion in \eqref{Cauchy2}
\begin{multline*}
\Rzeta(s)=\sum_{n=1}^\ell n^{-s}+\frac{(-1)^\ell}{2i}\xi^{-s}e^{\pi i \xi^2}\Bigl(
\sum_{k=0}^K \tau^k\int_{0\swarrow}\frac{e^{\frac{\pi i}{2}(v-q)^2}}
{2\cos\frac{\pi v}{2}}P_k\bigl(i\sqrt{\pi}(v-q)\bigr)\,dv\Bigr.\\
+\int_{0\swarrow}\frac{e^{\frac{\pi i}{2}(v-q)^2}}
{2\cos\frac{\pi v}{2}}Rg_K\bigl(\tau,i\sqrt{\pi}(v-q)\bigr)\,dv\Bigl.\Bigr).
\end{multline*}
By the definition of $\tau$ we may write this as 
\begin{equation}\label{RSformula}
\Rzeta(s)=\sum_{n=1}^\ell n^{-s}+(-1)^\ell \frac{1}{2i}\xi^{-s} 
e^{\pi i \xi^2}\Bigl\{\sum_{k=0}^K\frac{D_k(q)}{\xi^k}+RS_K\Bigr\}
\end{equation}
where
\begin{gather}\label{DefD}
D_k(q):=\Bigl(-\frac{1}{4\sqrt{\pi}}\Bigr)^k\int_{0\swarrow}
\frac{e^{\frac{\pi i}{2}(v-q)^2}}{2\cos\frac{\pi v}{2}}P_k
\bigl(i\sqrt{\pi}(v-q)\bigr)\,dv\\
RS_K:=\int_{0\swarrow}
\frac{e^{\frac{\pi i}{2}(v-q)^2}}{2\cos\frac{\pi v}{2}}Rg_K
\bigl(\tau,i\sqrt{\pi}(v-q)\bigr)\,dv\label{error}
\end{gather}

\section{Formulas for the computation of the terms.}\label{sectres}
The first $D_0(q)$ may be computed explicitly. (This is due to Riemann \cite{Siegel},
and a proof may be found in \cite{Chandra})
\begin{equation}\label{defvalueF}
D_0(q)=\int_{0\swarrow}\frac{e^{\frac{\pi i}{2}(v-q)^2}}
{2\cos\frac{\pi v}{2}}\,dv=\frac{e^{\frac{\pi i}{2}q^2}-\sqrt{2}\,e^{\frac{\pi i}{8}}\cos\frac{\pi q}{2}
}{\cos\pi q}=:G(q).
\end{equation}
We may compute the derivative of $G$ (with respect to $q$) differentiating under 
the integral sign. From the definition of the Hermite polynomials,
\begin{equation}
D^k(e^{-x^2})=(-1)^k e^{-x^2}H_k(x)
\end{equation}
we get (here $D$ denotes differentiation with respect to $q$)
\begin{equation}
D^k\bigl(e^{\frac{\pi i }{2}(v-q)^2}\bigr)=e^{\frac{\pi i }{2}(v-q)^2}
\Bigl(\frac{e^{-\frac{\pi i}{4}}\sqrt{\pi}}{\sqrt{2}}\Bigr)^k H_k\Bigl(
\frac{e^{-\frac{\pi i}{4}}\sqrt{\pi}}{\sqrt{2}}(v-q)\Bigr)
\end{equation}
It follows that
\begin{equation}
G^{(k)}(q)=\Bigl(\frac{e^{-\frac{\pi i}{4}}\sqrt{\pi}}{\sqrt{2}}\Bigr)^k
\int_{0\swarrow}\frac{e^{\frac{\pi i}{2}(v-p)^2}}{2\cos\frac{\pi v}{2}}
H_k\Bigl(\frac{e^{-\frac{\pi i}{4}}\sqrt{\pi}}{\sqrt{2}}(v-q)\Bigr)\,dv.
\end{equation}

From \eqref{defg} and \eqref{gexpansion} it is easy to see that the polynomial
$P_k(z)$ has degree $3k$, and it is also an even or an odd polynomial. We may expand $P_k(z)$ in terms of Hermite polynomials. 
Thus, there are certain numbers $d^{(k)}_j$ such that
\begin{equation}\label{intrdkj}
\frac{\pi^{2k}P_k(i\sqrt{\pi}(v-q))}{(-4\sqrt{\pi})^k}
=\sum_{j=0}^{\lfloor3k/2\rfloor}\Bigl(\frac{\pi }{2i}\Bigr)^j
d^{(k)}_j\Bigl(\frac{e^{-\frac{\pi i}{4}}\sqrt{\pi}}{\sqrt{2}}\Bigr)^{3k-2j}
H_{3k-2j}\Bigl(\frac{e^{-\frac{\pi i}{4}}\sqrt{\pi}}{\sqrt{2}}(v-q)\Bigr).
\end{equation}

From this we get the formula for $D_k(q)$
\begin{equation}\label{coefficients}
D_k(q)=\frac{1}{\pi^{2k}}\sum_{j=0}^{\lfloor3k/2\rfloor}\Bigl(\frac{\pi }{2i}\Bigr)^j
d^{(k)}_jG^{(3k-2j)}(q).
\end{equation}

Now we show how to compute the numbers $d^{(k)}_j$. First, observe that 
\eqref{intrdkj} may be written as 
\begin{equation}\label{defUk}
U_k(x):=\Bigl(\frac{e^{-\frac{\pi i}{4}}}{\sqrt{2}}\Bigr)^kP_k(x\sqrt{2}
e^{\frac{3\pi i}{4}})=\sum_{j=0}^{\lfloor3k/2\rfloor}d^{(k)}_jH_{3k-2j}(x).
\end{equation}
These polynomials $U_k(x)$ satisfy the following
\begin{equation}
\sum_{k=0}^\infty U_k(x)\tau^k=g\Bigl(\frac{e^{-\frac{\pi i}{4}}}{\sqrt{2}}\tau,
x\sqrt{2}e^{\frac{3\pi i}{4}}\Bigr)=\exp\Bigl\{\frac{1}{4\tau^2}\log(1-2x\tau)
+\frac{x}{2\tau}+\frac{x^2}{2}\Bigr\}
\end{equation}
To determine the polynomials $U_k(x)$, we differentiate logarithmically with
respect to $x$ and get
\begin{equation}
\sum_{k=0}^\infty U'_k(x)\tau^k=-\frac{2x^2\tau}{1-2x\tau}\sum_{k=0}^\infty
U_k(x)\tau^k
\end{equation}
from which we derive the relation
\begin{equation}
U'_k(x)=-2x^2U_{k-1}(x)+2xU'_{k-1}(x).
\end{equation}
We extend the definition of $d^{(k)}_j$ so that $d^{(k)}_j=0$ for $j<0$ and 
for $j>3k/2$. Then by \eqref{defUk}  we have 
\begin{equation}
\sum_j  d^{(k)}_jH'_{3k-2j}(x)=-2x^2\sum_jd^{(k-1)}_j H_{3k-3-2j}(x)
+2x\sum_jd^{(k-1)}_jH'_{3k-3-2j}(x).
\end{equation}
It is well known that the Hermite polynomials satisfy
\begin{equation}
H'_n(x)=2nH_{n-1}(x),\quad xH_n(x)=\tfrac12  H_{n+1}(x)+nH_{n-1}(x)
\end{equation}
so that 
\begin{equation}
x^2H_n(x)=\tfrac14  H_{n+2}(x)+(n+\tfrac12 )H_n(x)+n(n-1)H_{n-2}(x).
\end{equation}
Therefore, the above relation can be written as 
\begin{multline}
\sum_j d^{(k)}_j(6k-4j)H_{3k-2j-1}=\\ 
-2\sum_jd^{(k-1)}_j\bigl\{\tfrac14  H_{3k-1-2j}+(3k-3-2j+\tfrac12 )H_{3k-3-2j}+
+(3k-3-2j)(3k-4-2j)H_{3k-5-2j}\bigr\}\\ 
+2\sum_jd^{(k-1)}_j2(3k-3-2j)\bigl\{\tfrac12  H_{3k-3-2j}+(3k-4-2j)H_{3k-5-2j}\bigr\}.
\end{multline}
Separating the coefficients of $H_{3k-2j-1}$ and simplifying, we arrive at the equation
\begin{equation}\label{coef_d1}
(6k-4j)d^{(k)}_j=-\tfrac12  d^{(k-1)}_j-d^{(k-1)}_{j-1}+2(3k-2j)(3k-2j+1)d^{(k-1)}_{j-2}.
\end{equation}
Starting from $d^{(0)}_0=1$, $d^{(0)}_j=0$ for $j\ne0$, these equations determine 
the $d^{(k)}_j$ except those with $3k=2j$. To determine these for $k>0$ we use  the
relation  \eqref{defUk} for $x=0$, where $P_k(0)=0$:
\begin{equation}\label{coef_d2}
d^{(k)}_{3k/2}=-\sum_{j=0}^{3k/2-1}(-1)^{3k/2-j}d^{(k)}_j\frac{(3k-2j)!}{(3k/2-j)!},
\qquad 3k\equiv0 \pmod{2}.
\end{equation}

\begin{remark}
The numbers $d^{(k)}_j$ are related to the ones defined in \cite{A86}. The numbers defined in \cite{A86} are polynomials in $\sigma$, denoting these polynomials by $d^{(k)}_j(\sigma)$ we have the relation $d^{(k)}_j=(-1)^k d^{(k)}_j(0)$. This can be  easily proved by comparing the recurrence definitions given here and in \cite{A86}.
\end{remark}

\begin{remark} 
Equation \eqref{coefficients} gives us the following explicit formulas:
\begingroup\makeatletter\def\f@size{9}\check@mathfonts
\def\maketag@@@#1{\hbox{\m@th\large\normalfont#1}}%
\begin{equation}\label{explicit}
\begin{aligned}
D_0(q)&=G(q),\\
D_1(q)&=\frac14\Bigl(\frac{i}{\pi}G'(q)-\frac{1}{3\pi^2}G^{(3)}(q)\Bigr),\\
D_2(q)&=\frac{1}{16}\Bigl(-\frac{2i}{3\pi}G(q)+\frac{1}{2\pi^2}G''(q)-\frac{i}{3\pi^3}G^{(4)}(q)+\frac{1}{18\pi^4}G^{(6)}(q)\Bigr),\\
D_3(q)&=\frac{1}{64}\Bigl(-\frac{1}{3\pi^2}G'(q)+\frac{31i}{18\pi^3}G^{(3)}(q)-\frac{11}{30\pi^4}G^{(5)}(q)+\frac{i}{18\pi^5}G^{(7)}(q)-\frac{1}{162\pi^6}G^{(9)}(q)\Bigr).
\end{aligned}
\end{equation}
\endgroup
\end{remark}

\section{Bound for the error.}

Next, we obtain an integral representation of the remainder $Rg_K(\tau,z)$ in the expansion 
\eqref{gexpansion}.

\begin{lemma}\label{L:Integral}
Let $0<r<1$ and $0<\phi<\pi$, denote by $L$ the line $r+xe^{i\phi}$ with $x\in\R$. Then for any any complex numbers $\tau\ne0$ and  $z$  such that 
$-2iz\tau$ is in the half plane limited by $L$ containing $0$ and any non negative integer $K$ we have  
\begin{equation}\label{RgK}
Rg_K(\tau,z)=\frac{(-2iz\tau)^{K+1}}{2\pi i}\int_L\frac{e^{-\frac{i z^2}{2}f(\zeta)}}
{(\zeta+2iz\tau)\zeta^K}\frac{d\zeta}{\zeta},
\end{equation}
where 
\begin{equation}\label{deff}
f(\zeta):=-\frac{\log(1-\zeta)}{\zeta^2}-\frac{1}{\zeta}-\frac{1}{2}
\end{equation}
\end{lemma}

\begin{proof}
For $z=0$, the result is trivial. Assuming that $z\ne0$, 
by \eqref{gexpansion} and Cauchy's Theorem we have
\begin{displaymath}
Rg_K(\tau,z)=\frac{1}{2\pi i}\int_C\frac{g(\zeta,z)}{\zeta-\tau}\Bigl(\frac{\tau}{\zeta}
\Bigr)^{K+1}\,d\zeta
\end{displaymath}
where 
$C$ is a simple path in the region where $\zeta\mapsto g(\zeta,z)$ is holomorphic and has $0$ and $\tau$ in its interior.
We change the variables by putting $i\zeta/2z$ instead of $\zeta$. In this way, we show
that
\begin{equation}
Rg_K(\tau,z)=\frac{(-2iz\tau)^{K+1}}{2\pi i}\int_C\frac{e^{-\frac{i z^2}{2}f(\zeta)}}
{(\zeta+2iz\tau)\zeta^K}\frac{d\zeta}{\zeta}
\end{equation}
where $f$ is  defined in \eqref{deff}, denoting by $\log$ the main branch of the logarithm
and now $C$ denotes a simple path that circles the points $0$ and $-2iz\tau$, 
and is contained in the region obtained from $\C$ by making a cut along the positive 
real axis from $1$ to $+\infty$ (the domain where $f$ is holomorphic). 

Note that $f(\zeta)$ is holomorphic on $\C\smallsetminus[1,+\infty)$ and is bounded in this domain except for a small neighborhood of $\zeta=1$. Therefore, a standard application of Cauchy's theorem  yields 
\[Rg_K(\tau,z)=\frac{(-2iz\tau)^{K+1}}{2\pi i}\int_L\frac{e^{-\frac{i z^2}{2}f(\zeta)}}
{(\zeta+2iz\tau)\zeta^K}\frac{d\zeta}{\zeta}.\qedhere\]
\end{proof}

\begin{definition}
For $q\in\C$ define 
\begin{equation}\label{defmunu}
\mu=\mu(q):=\frac{\Re(q)+\Im(q)}{\sqrt{2}},\quad \nu=\nu(q)=\frac{\Im(q)-\Re(q)}{\sqrt{2}}.
\end{equation}
So that $q=(\mu+i\nu)e^{\pi i/4}$.
\end{definition}

\begin{theorem}\label{Theorem1}
Let $0<\theta<\pi$ be given, and let 
\begin{equation}\label{E:Delta}
\Delta=\{s\in\C: -\pi+\theta\le\arg(s)\le\pi-\theta, |s|\ge2\pi\}.
\end{equation}
For $s\in\Delta$ and $K\ge1$ a natural number,  define $\xi$ as in \eqref{defeta},  $\ell$ as in \eqref{defm}, $q$ as in \eqref{defp}, $\mu$ as in \eqref{defmunu}, $D_k(q)$ as in \eqref{DefD} and $RS_K$ as in \eqref{error}.
For  $s\to\infty$ we have
\begin{equation}\label{firstexpansion}
\Rzeta(s)=\sum_{n=1}^\ell n^{-s}+\frac{(-1)^\ell}{2i}\xi^{-s} 
e^{\pi i \xi^2}\Bigl\{\sum_{k=0}^K\frac{D_k(q)}{\xi^k}+
\Orden_{K,\theta}(e^{-\frac{\pi}{2\sqrt{2}}|\mu|}|\xi|^{-K-1})\Bigr\}.
\end{equation}
\end{theorem}

\begin{proof}
Since we restrict $s$ to $\Delta$  we  have $|\xi|\ge1$ and $-\frac{3\pi}{4}+\frac{\theta}{2}\le\arg(
\xi)\le\frac{\pi}{4}-\frac{\theta}{2}$. We have defined $\ell$ as the only integer with
$\ell\le\Re\xi-\Im\xi<\ell+1$. It follows that $\ell\ge0$. Let us define $q_1$ and 
$q_2$ as the real and imaginary parts of $q=q_1+iq_2$. Since $q=-2\ell-1+2\xi$
we have (putting $\xi=\xi_1+i\xi_2$) 
\[ q_1-q_2=2\xi_1-2\ell-1-2\xi_2=2(\xi_1-\xi_2-\lfloor\xi_1-\xi_2\rfloor-\tfrac12)\in[-1,1)\]
So, $q$ is contained in the strip $-1\le q_1-q_2<1$.  

Recall that the error is given by \eqref{error}
\[RS_K=\int_{0\swarrow}
\frac{e^{\frac{\pi i}{2}(v-q)^2}}{2\cos\frac{\pi v}{2}}Rg_K
\bigl(\tau,i\sqrt{\pi}(v-q)\bigr)\,dv\]
where the path is a line through the origin  of the direction $e^{-3\pi i/4}$.
By Cauchy's theorem we may change the line of integration to any parallel if we do not pass through 
the singularity of $Rg_K$ at the point $v$ where $2\tau \sqrt{\pi}(v-q)=1$ nor through the poles at $v=-1$ and $v=1$. 

Pick a small and positive $\beta$ (we shall see that a convenient election is 
$\beta=\frac{1}{3}\sin\frac{\theta}{2}$). Then the path of integration $0\swarrow$ will be changed to a parallel line $R$ in the following way
If $-1+{\beta}{\sqrt{2}}\le q_1-q_2\le 1-{\beta}{\sqrt{2}}$ we take $v=q+xe^{-3\pi i/4}$. In this case, the path will pass exactly through the point $q$. 
If  $-1\le q_1-q_2<-1+{\beta}{\sqrt{2}}$ then we take  
$v=q+xe^{-3\pi i/4}+\beta e^{-\pi i/4}$.
If $1-{\beta}{\sqrt{2}}<q_1-q_2<1$ then we take 
$v=q+xe^{-3\pi i/4}-\beta e^{-\pi i/4}$. 
So in all cases we will have a path $R$ given by 
$v=q+xe^{-3\pi i/4}+\gamma e^{-\pi i/4}$,  where $\gamma$ may be $\pm\beta$ or $0$.  

For the points $v$ of this line $R$ we have in the integrand $Rg_K(\tau,z)$ where $z=i\sqrt{\pi}(v-q)$. Therefore,  for 
$z= i\sqrt{\pi}(xe^{-3\pi i/4}+\gamma e^{-\pi i/4})$ with $x\in \R$. For these values, 
we have
\[-2iz\tau=2\sqrt{\pi}(v-q)\tau=2\sqrt{\pi}(xe^{-3\pi i/4}+\gamma e^{-\pi i/4}) \frac{-1}{4\sqrt{\pi}\xi}=\frac{xe^{\pi i/4}}{2\xi}+\frac{\gamma e^{3\pi i/4}}{2\xi}=\frac{x+i\gamma}{2|\xi|}e^{i\phi},
\]
where $\phi=\frac{\pi}{4}-\arg\xi$, is contained in $[\frac{\theta}{2},\pi-\frac{\theta}{2}]$.
Hence $-2iz\tau$ is contained in a line of direction $e^{i\phi}$ that cuts the real axis at the point $-\frac{\gamma}{2|\xi|\sin\phi}$. The absolute value of this is
\[\Bigl|-\frac{\gamma}{2|\xi|\sin\phi}\Bigr|=\frac{\beta}{2|\xi|\sin\phi}\le 
\frac{\beta}{2|\xi|\sin\frac{\theta}{2}}=\frac{1}{6|\xi|}\le \frac16,\]
since $|\xi|\ge 1$.   In this way the line $R$ satisfies the needed conditions to apply Cauchy's Theorem so that 
\begin{equation}\label{error_modified}
RS_K:=\int_{R}
\frac{e^{\frac{\pi i}{2}(v-q)^2}}{2\cos\frac{\pi v}{2}}Rg_K
\bigl(\tau,i\sqrt{\pi}(v-q)\bigr)\,dv
\end{equation}

Hence, taking $L$ as a line with direction $e^{i\phi}$ and passing through the point $1/2$, Lemma \ref{L:Integral} applies to any value of $z$ that interests us.

Substituting \eqref{RgK} into \eqref{error_modified} and applying Fubini's Theorem, we get
\begin{equation}\label{doubleInt}
RS_K=\frac{(2\tau\sqrt{\pi})^{K+1}}{2\pi i}\int_L\frac{d\zeta}{\zeta^{K+1}}
\int_{R}\frac{e^{\frac{\pi i}{2}(v-q)^2}}{2\cos\frac{\pi v}{2}}
\frac{(v-q)^{K+1}e^{\frac{\pi i}{2}(v-q)^2f(\zeta)}}{
(\zeta-2\sqrt{\pi}(v-q)\tau)}\,dv
\end{equation}
First, we shall bound the inner integral. 

We bound the terms of the integrand. Since $\zeta\in L$ we have $\zeta=1/2+y e^{i\phi}$ with $y\in\R$ and
\[
|\zeta-2\sqrt{\pi}(v-q)\tau|=\Bigl|\frac12+ye^{i\phi}-\frac{x+i\gamma}{2|\xi|}e^{i\phi}\Bigr|
=\Bigl|\frac12+\Bigl(y-\frac{x}{2|\xi|}\Bigr)e^{i\phi}-\frac{\gamma i }{2|\xi|}
e^{i\phi}\Bigr|
\]
This is greater than
\[\ge \min_{x\in\R}\Bigl|\frac12+xe^{i\phi}\Bigr|-\frac{\beta}{2|\xi|}
\ge\frac{1}{2}\sin\phi-\frac{\beta}{2}\ge
\frac{1}{2}\sin\frac{\theta}{2}-\frac{1}{6}\sin\frac{\theta}{2}=
\frac{1}{3}\sin\frac{\theta}{2}.
\]
 Therefore,
\begin{equation}
|\zeta-2\sqrt{\pi}(v-q)\tau|\ge\frac{1}{3}\sin\frac{\theta}{2}
\end{equation}

The point $v$ of our path of integration is always on the strip 
$-1+{\beta}{\sqrt{2}}\le v_1-v_2\le 1-{\beta}{\sqrt{2}}$. 
Since $\cos\frac{\pi v}{2}\to\infty$ when $v\to\infty$ in this strip, and it is 
a continuous function that does not vanish on this strip, we will get a 
constant $M(\theta)>0$ such that 
\begin{equation}
\Bigl|\cos\frac{\pi v}{2}\Bigr|\ge M(\theta), 
\qquad -1+{\beta}{\sqrt{2}}\le v_1-v_2\le 1-{\beta}{\sqrt{2}}.
\end{equation}

This can be improved for $\mu$ large. In fact,
\begin{align*}
\cos\frac{\pi v}{2}&=\cos\frac{\pi}{2}\bigl[(\mu+i\nu)e^{\frac{\pi i}{4}}-(x+i\gamma)e^{\frac{\pi i}{4}}\bigr]\\
&=\frac{\pi}{2\sqrt{2}}\bigl[(\mu-x)-(\nu-\gamma)+i((\mu-x)+(\nu-\gamma))\bigr]
\end{align*}
For real numbers $a$ and $b$,  $e^{-|b|}|\cos(a+ib)|\ge \frac12(1-e^{-2|b|})$. Hence, 
\[\Bigl|\cos\frac{\pi v}{2}\Bigr|\ge \frac14\exp\Bigl(\frac{\pi}{2\sqrt{2}}|(\mu-x)+(\nu-\gamma)|\Bigr)\]
for $|(\mu-x)+(\nu-\gamma)|>1$.  Since $|\gamma|\le|\beta|\le \frac13$ and $|\nu|=\frac{1}{\sqrt{2}}|q_1-q_2|\le \frac{1}{\sqrt{2}}$, there is an absolute constant $c_0$ with
\[\Bigl|\cos\frac{\pi v}{2}\Bigr|\ge c_0 e^{\frac{\pi}{2\sqrt{2}}|\mu-x|},\qquad |\mu-\nu|>c_0.\]
Combining the two inequalities, there is some constant $M(\theta)$ (perhaps different from the previous one such that 
\begin{equation}
\Bigl|\cos\frac{\pi v}{2}\Bigr|\ge M(\theta) e^{\frac{\pi}{2\sqrt{2}}|\mu-x|},\qquad v\in R.
\end{equation}

\begin{displaymath}
\frac{\pi i}{2}(v-q)^2=\frac{\pi i}{2}(xe^{-\frac{3\pi i}{4}}+
\gamma e^{-\frac{\pi i}{4}})^2=-\frac{\pi x^2}{2}+\frac{\pi\gamma^2}{2}
-\pi\gamma ix.
\end{displaymath}
Therefore,
\begin{equation}
|e^{\frac{\pi i}{2}(v-q)^2}|\le e^{\frac{\pi\beta^2}{2}}e^{-\frac{\pi x^2}{2}}.
\end{equation}
The point $\zeta-\frac12=ye^{i\phi}$,  so it is in one of the angles
$-\pi+\frac{\theta}{2}<\arg(\zeta-\tfrac12 )<-\frac{\theta}{2}$ or $\frac{\theta}{2}<\arg(\zeta-\tfrac12 )<
\pi-\frac{\theta}{2}$.  Put $f(\zeta)=V(\zeta)+iU(\zeta)$. The two harmonic functions
are continuous on the closed angle and $\lim_{\zeta\to\infty}f(\zeta)=-1/2$. Therefore,  there are constants $V_0=V_0(\theta)$ and $U_0=U_0(\theta)$ (only depending on $\theta$) such that
for every $\zeta$ we  have $|U(\zeta)|\le U_0$ and  $|V(\zeta)|\le V_0$. These values are taken at the boundary of the angles. 

Then we have 
\begin{displaymath}
|e^{\frac{\pi i}{2}(v-q)^2 f(\zeta)}|= e^{-\frac{\pi x^2}{2}V}
e^{\frac{\pi \gamma^2}{2}V}e^{\pi \gamma U x}\le e^{\frac{\pi \beta^2}{2}V_0}
e^{-\frac{\pi x^2}{2}V}e^{\pi\beta U_0|x|}
\end{displaymath}
Therefore, the inner integral in \eqref{doubleInt} is bounded as 
\begin{equation}
\Bigl|\int_{R}\cdots\Bigr|\le \frac{4 e^{\frac{\pi\beta^2}{2}(1+V_0)}
}{M(\theta)\sin\frac{\theta}{2}}\int_{-\infty}^{+\infty}e^{-\frac{\pi x^2}{2}(1+V)}
e^{\pi \beta U_0 |x|-\frac{\pi}{2\sqrt{2}}|\mu-x|} (x^2+\beta^2)^{\frac{K+1}{2}}\,dx
\end{equation}
It is important here to note that $V>-1$ for all values of $\zeta$. This  
easily follows from the fact that $V(\overline\zeta)=V(\zeta)$ and  $V$ is a harmonic function on the upper half-plane, $\lim_{\zeta\to\infty} V(\zeta)=-\frac12$,  and on the real axis it is always $V(x)\ge-1$ the equality only in the case $x=2$. Thus, the above integral is finite and  there is a function $a(\theta)>0$ with 
$a(\theta)<1+V(\zeta)$ for all the values of $\zeta$ that interest us. 

Applying Lemma \ref{L:simplebound} in the appendix, we find a constant $C(\theta)$ only depending on $\theta$
such that
\begin{displaymath}
\int_{-\infty}^{+\infty} e^{-\frac{\pi x^2}{2}a(\theta)}e^{2\pi \beta U_0 |x|-\frac{\pi}{\sqrt{2}}|\mu-x|}\,dx \le  C(\theta)^2e^{-\frac{\pi}{\sqrt{2}}|\mu|},
\end{displaymath}
and by Schwarz inequality we get 
\begin{displaymath}
\Bigl|\int_{R}\cdots\Bigr|\le
C(\theta)e^{-\frac{\pi}{2\sqrt{2}}|\mu|}\Bigl(\int_{-\infty}^{+\infty}(x^2+\beta^2)^{K+1}e^{-\frac{\pi x^2}{2}a(\theta)}\,dx\Bigr)^{1/2}
\end{displaymath}
changing here $(x^2+\beta^2)^{K+1}\le (2x^2)^{K+1}+(2\beta^2)^{K+1}$ we get 
\begin{equation}
\Bigl|\int_{R}\cdots\Bigr|\le
C(\theta)\Bigl(\frac{(2\beta^2)^{K+1}}{\sqrt{2a(\theta)}}+\frac{2^K\Gamma(\frac32+K)}
{(\frac{\pi}{2}a(\theta))^{\frac32+K}}\Bigr)^{1/2}e^{-\frac{\pi}{2\sqrt{2}}|\mu|}
\end{equation}
Let us call $C_K(\theta)$ this bound. 
Then
\begin{displaymath}
|RS_K|\le C_K(\theta)e^{-\frac{\pi}{2\sqrt{2}}|\mu|}\frac{(2\tau\sqrt{\pi})^{K+1}}{2\pi}\int_L\frac{|d\zeta|}{|\zeta|^{K+1}}
\end{displaymath}
This integral is bounded when $K\ge1$.  So, by the definition of $\tau$ we get
\begin{equation}
|RS_K|\le \frac{C_K(\theta)e^{-\frac{\pi}{2\sqrt{2}}|\mu|}}{|\xi|^{K+1}}
\end{equation}
with a new constant $C_K(\theta)$. 
\end{proof}

\section{Properties of the function \texorpdfstring{$G(q)$}{G(q)}.}
To use the asymptotic expansion, we need to consider the function $G(q)$. Recall that 
we are interested in the values of $q$ contained in the strip $-1\le \Re(q)-\Im(q)\le 1$.  Here we will consider  more general strips. For $a>0$ we define 
\begin{equation}\label{defBa}
B_a=\{(\mu+i\nu)e^{\pi i/4}\colon \mu, \nu\in\R, |\nu|\le a/\sqrt{2}\}
\end{equation}
Hence $q=q_1+iq_2\in B_1$ is equivalent to $|q_1-q_2|\le 1$. 

Let $B_a^+$ denote the set of $q=(\mu+i\nu)e^{\pi i/4}\in B_a$ with $\mu\ge0$. 
Since $(\mu+i\nu)e^{\pi i/4}\in B_a$ is equivalent to $(-\mu-i\nu)e^{\pi i/4}\in B_a$
and $G(q)$ is an even function, the behavior of $G(q)$ in $B_a$ is determined by its behaviour on $B_a^+$.

\begin{proposition}\label{Gzeros}
The function $G(q)$ does not vanish on the strip $B_1$. 
\end{proposition}

\begin{proof}
Assume that $q\in B_1$, by definition, there exist real numbers $\mu$ and $\nu$ with
$q=(\mu+i\nu)e^{\pi i/4}$ and $|\nu|\le1/\sqrt{2}$. First, we show that the function does not vanish for $|\mu|\ge2$. Since $G(-q)=G(q)$ we may assume that $\mu>2$. In this case, we have
\[
|e^{\frac{\pi i}{2}q^2}|=|e^{-\frac{\pi}{2}(\mu+i\nu)^2}|=e^{-\frac{\pi}{2}(\mu^2-\nu^2)}\le e^{\frac{\pi}{4}}e^{-\frac{\pi}{2} \mu^2},\]
\[|e^{\frac{\pi i}{2} q}|=|e^{\frac{\pi}{2}(\mu+i\nu)(\frac{-1+i}{\sqrt{2}})}|=
e^{\frac{\pi}{2\sqrt{2}}(-\mu-\nu)}\le e^{\frac{\pi}{4}}e^{-\frac{\pi \mu}{2\sqrt{2}}},\]
\[|e^{-\frac{\pi i}{2} q}|=e^{\frac{\pi}{2\sqrt{2}}(\mu+\nu)}\ge e^{-\frac{\pi}{4}}e^{\frac{\pi \mu}{2\sqrt{2}}}.\]

It follows that
\begin{align*}
|e^{\frac{\pi i}{2}q^2}-\sqrt{2}e^{\frac{\pi i}{8}}\cos\tfrac{\pi}{2}q|&\ge
\frac{1}{\sqrt{2}}e^{-\frac{\pi}{4}}e^{\frac{\pi \mu}{2\sqrt{2}}}-
\frac{1}{\sqrt{2}}e^{\frac{\pi}{4}}e^{-\frac{\pi \mu}{2\sqrt{2}}}-e^{\frac{\pi}{4}}e^{-\frac{\pi}{2} \mu^2}\\
&=\frac{1}{\sqrt{2}}\sinh\Bigl(\frac{\pi \mu}{2\sqrt{2}}-\frac{\pi}{4}\Bigr)-e^{\frac{\pi}{4}}e^{-\frac{\pi}{2} \mu^2}
\end{align*}
This is $>0$ for $\mu=2$, and since it is increasing for all $\mu\ge2$. 

\begin{figure}[H]
\begin{center}
\includegraphics[width=\hsize]{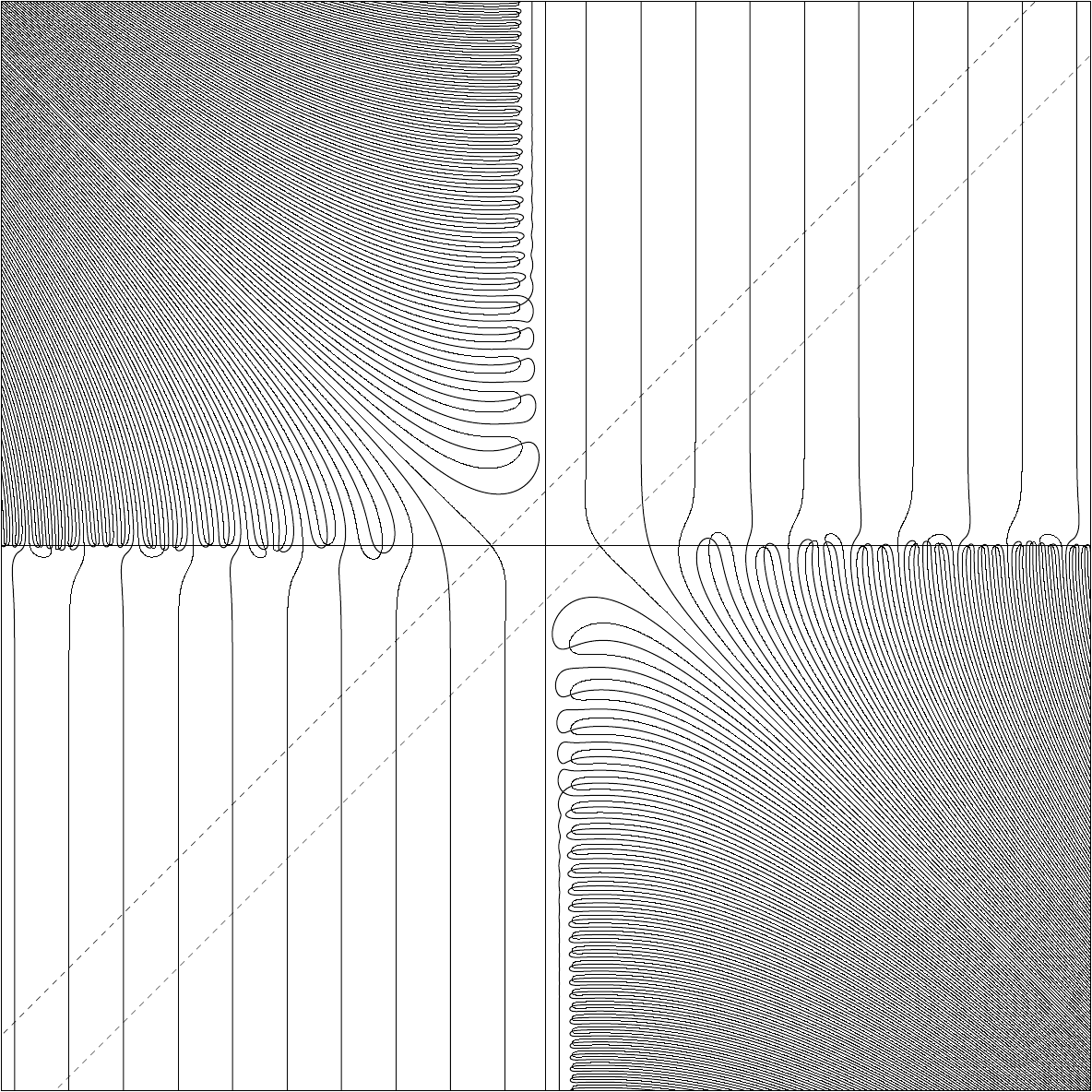}
\caption{x-ray of $G(q)$ for $q\in(-10,10)^2$ with  the region $B_1$ marked. }
\label{functG}
\end{center}
\end{figure}

The image of the parallelogram determined by $|\mu|\le 2$ and $|\nu|\le 1/\sqrt{2}$ is given in the figure \ref{border}.
Since the image of the border does not enclose $0$ the function does not take the value $0$ on this parallelogram. 
This is also shown in the x-ray of the function $G(q)$. Both can be transformed into a conventional computational proof. 
\end{proof}

\begin{proposition}\label{Gasymptotic}
For  $q=(\mu+i\nu)e^{\pi i/4}\in B_a^+$,  and $q\to\infty$  we have
\begin{equation}\label{G-approx}
G(q)=-\sqrt{2}e^{\frac{\pi i}{8}} \exp\Bigl(-\frac{\pi}{2\sqrt{2}}(\mu+\nu-i(\mu-\nu)\Bigr)(1+\Orden_a(e^{-\pi \mu/\sqrt{2}})).
\end{equation}
\end{proposition}

\begin{proof}
Let  $q=(\mu+i\nu)e^{\pi i/4}\in B_a$ with $\mu>0$ be given, then  $|e^{\frac{\pi i}{2}q^2}|=\exp(-\frac{\pi \mu^2}{2}+\frac{\pi \nu^2}{2})$,
and we will see that this term is of lesser order than $\cos\frac{\pi q}{2}$ when $q\in B_a$ and $\mu$ large.

We have
\begin{align*}
\cos\frac{\pi q}{2}&=\frac12\exp\Bigl(\frac{\pi}{2\sqrt{2}}(-\mu-\nu+i(\mu-\nu)\Bigr)+
\frac12\exp\Bigl(\frac{\pi}{2\sqrt{2}}(\mu+\nu-i(\mu-\nu)\Bigr)\\
&=\frac12\exp\Bigl(\frac{\pi}{2\sqrt{2}}(\mu+\nu-i(\mu-\nu)\Bigr)(1+\Orden_a(e^{-\pi \mu/\sqrt{2}}))
\end{align*}

\begin{figure}[H]
\begin{center}
\includegraphics[width=\hsize]{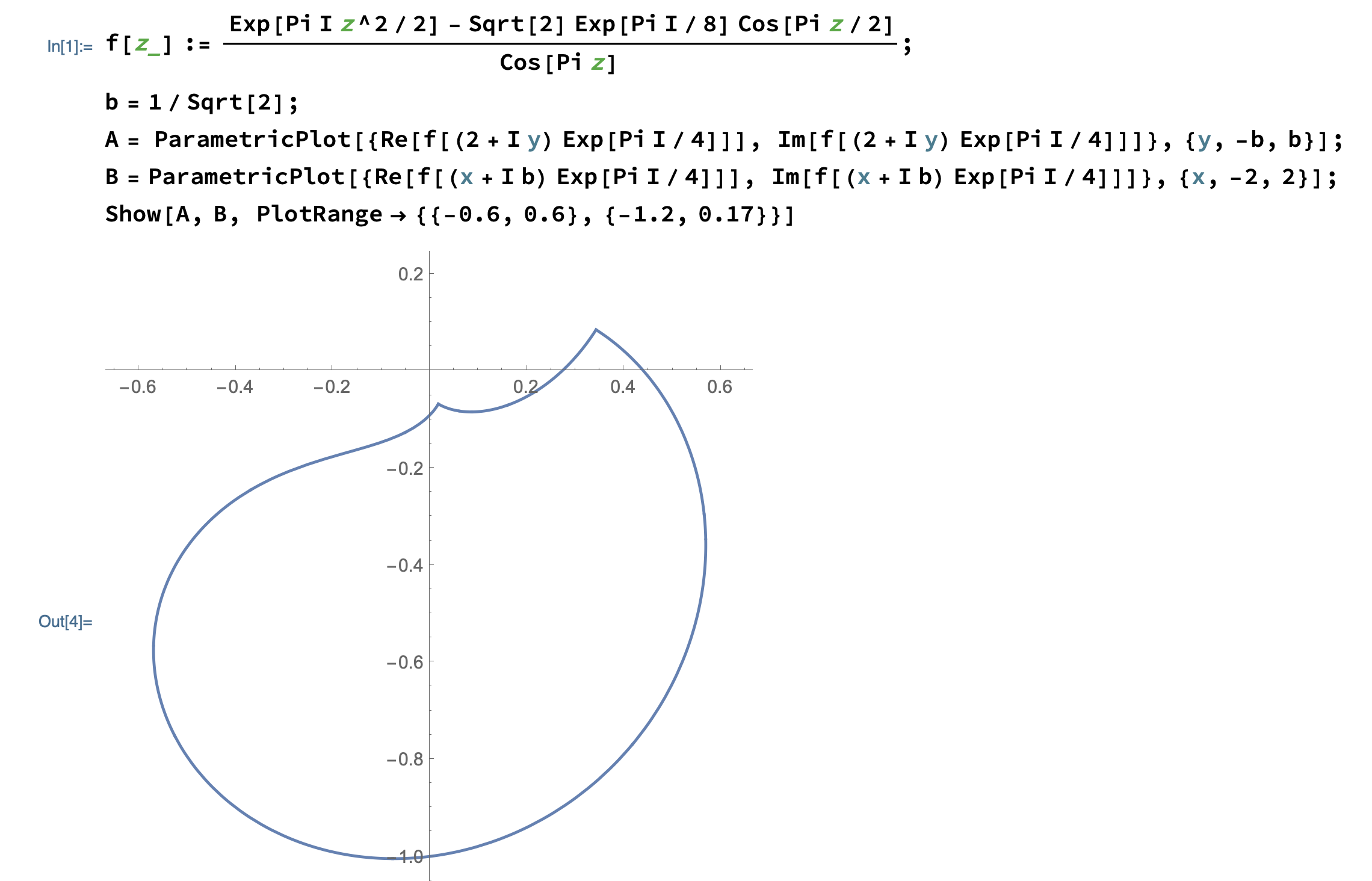}
\caption{image of the region $|\mu|\le 2$, $|\nu|\le 1/\sqrt{2}$.}
\label{border}
\end{center}
\end{figure}

Analogously,
\[\cos\pi q=\frac12\exp\Bigl(\frac{\pi}{\sqrt{2}}(\mu+\nu-i(\mu-\nu)\Bigr)(1+\Orden_a(e^{-\pi \mu\sqrt{2}}))\]
Hence, the explicit form of $G(q)=\frac{e^{\frac{\pi i}{2}q^2}-\sqrt{2}\,e^{\frac{\pi i}{8}}\cos\frac{\pi q}{2}
}{\cos\pi q}$ yields \eqref{G-approx}.
\end{proof}

\begin{proposition}\label{Gproperties}
(a) For any given $a>0$, there exist constants $c_0(a)>0$ and $C_0(a)$ such that 
for $q=(\mu+i\nu)e^{\pi i/4}\in B_a$ we have
\begin{equation}\label{E:Gineq1}
c_0(a)e^{-\frac{\pi}{2\sqrt2}|\mu|}\le |G(q)|\le C_0(a)e^{-\frac{\pi}{2\sqrt2}|\mu|}
\end{equation}
(b) For any given $a>0$, and non-negative integer $k$ there exist constants $C_k(a)$ such that 
\begin{equation}\label{E:Gineq2}
|G_k(q)|\le C_k(a) e^{-\frac{\pi}{2\sqrt{2}}|\mu|},\qquad q\in B_a.
\end{equation}
\end{proposition}

\begin{proof}
As noted above, we only have to prove (a) and (b) for $x>0$.  By Proposition \ref{Gasymptotic}  the function $G(q)e^{\frac{\pi}{2\sqrt2}\mu}$ is bounded in $B_a^+$.  For $q\to\infty$ in $B_a^+$ we have by \eqref{G-approx}
\[G(q)e^{\frac{\pi}{2\sqrt2}\mu}=-\sqrt{2}e^{\frac{\pi i}{8}} \exp\Bigl(-\frac{\pi}{2\sqrt{2}}(\nu-i(\mu-\nu)\Bigr)(1+\Orden_a(e^{-\pi \mu/\sqrt{2}})).\]
For $q\in B_a^+$ we have $|\nu|\le a/\sqrt2$, therefore given $\varepsilon\in(0,1)$ there is some $\mu_0>0$ such that for $\mu>\mu_0$ 
\[(1-\varepsilon)\sqrt{2}e^{-\frac{\pi a}{4}}\le |G(q)|e^{\frac{\pi}{2\sqrt2}\mu}\le (1+\varepsilon)\sqrt{2}e^{\frac{\pi a}{4}}.\]
The continuous function $|G(q)e^{\frac{\pi}{2\sqrt2}\mu}|$ is bounded
in the compact set  $ B_a^+$ and $0\le \mu\le \mu_0$. These two observations prove the inequality (b) for $k=0$.  Inequality (b) for $k\ge1$ follows from Cauchy's formula
\[G^{(k)}(q)=\frac{k!}{2\pi i}\int_C\frac{G(\zeta)}{(\zeta-q)^{k+1}}\,d\zeta.\]
If $q\in B_a$, taking $C$ as a circle with center at $q$ and radius 1, we have $\zeta\in B_{a+1}$. Since $|\mu(\zeta)|\ge |\mu(q)|-1$ it follows that
\[|G^{(k)}(q)|\le k! C_0(a+1)e^{\frac{\pi}{2\sqrt{2}}}e^{-\frac{\pi}{2\sqrt{2}}|\mu|}.\]

It remains only to prove the left inequality in (a). By the above reasoning, we have proved that 
\[(1-\varepsilon)\sqrt{2}e^{-\frac{\pi a}{4}}\le |G(q)|e^{\frac{\pi}{2\sqrt2}\mu},\qquad q\in B_1^+, \mu>\mu_0.\]
In the compact set $B_1^+$, $0\le \mu\le \mu_0$, the continuous function $|G(q)|e^{\frac{\pi}{2\sqrt2}\mu}$ takes a minimum value that according to Proposition \ref{Gzeros} is strictly positive. This is all we need. \end{proof}

\section{Asymptotic behaviour on the right.}

The asymptotic expansion simplifies in a closed set $L$ containing the first  and most of the fourth quadrant.  The next Lemma is useful to define $L$.

\begin{lemma}\label{L:region}
Let $P\subset\C$ defined by 
\begin{equation}
P=\{\xi\in \C\colon |\xi|\ge e, \quad|\xi^{-2\pi i\xi^2}e^{\pi i\xi^2+\pi\xi}|\le 1,\quad -\tfrac{\pi}{2}\le \arg\xi\le -\tfrac{\pi}{4}\},
\end{equation}
where  $\xi^z=\exp(z\log\xi)$ is defined taking  $\log\xi=\log|\xi|+i\arg\xi$ and $-\pi/2\le \arg\xi\le 0$. In $P$ we have  $\xi=\xi_1+i\xi_2$, with $\xi_1\ge0$ and $\xi_2\le0$.

There is a function $\varphi\colon[e,+\infty)\to[0,\pi/4]$ such that 
\begin{equation}
P=\{\xi\in\C\colon|\xi|\ge e,\quad -\tfrac{\pi}{2}+\varphi(|\xi|)\le \arg\xi\le -\pi/4\}.
\end{equation}
\end{lemma}

\begin{proof}
Let $\xi=\xi_1+i\xi_2$ be in the fourth quadrant with $-\frac{\pi}{2}\le \arg\xi\le-\frac{\pi}{4}$, then 
\[\Re(-2\pi i\xi^2\log\xi+\pi i\xi^2+\pi\xi)=4\pi\xi_1\xi_2\log|\xi|+2\pi(\xi_1^2-\xi_2^2)\arg\xi-2\pi\xi_1\xi_2+\pi\xi_1.\]
Let $\varphi=\arctan(-\xi_1/\xi_2)$. Polar coordinates $(r,\varphi)$
satisfies
\[0\le \varphi\le \tfrac{\pi}{4},\quad\xi_1=r\sin\varphi,\quad \xi_2=-r\cos\varphi, \quad |\xi|=r\ge e,\quad \arg\xi=\varphi-\tfrac{\pi}{2}.\]
In terms of these coordinates, a point $\xi\in P$ just when 
$r\ge e$ and 
\[u(r,\varphi):=(2\log r)\sin(2\varphi)-2(\tfrac{\pi}{2}-\varphi)\cos(2\varphi)-\sin(2\varphi)-\frac{\sin\varphi}{r}\ge 0.\]
For $r\ge e$ we have $u(r,0)=-\pi<0$ and $u(r,\pi/4)=2\log r-1-\frac{1}{\sqrt{2} r}>0$.  Therefore, there is a value $\varphi(r)$ such that $u(r,\varphi(r))=0$. It can be shown that $\frac{du}{d\varphi} >0$ for $0\le\varphi\le\frac{\pi}{4}$,  hence $\varphi(r)$  is unique and $u(r,\varphi)\ge0$  just for $-\frac{\pi}{4}\ge \arg\xi\ge\varphi(r)-\frac{\pi}{2}$.
\end{proof}

\begin{proposition}\label{orderphi}
For $r\to+\infty$, the function $\varphi(r)$ defined in Lemma \ref{L:region} satisfies 
\begin{equation}
\lim_{r\to\infty}\frac{4\log r}{\pi} \sin\varphi(r)=1.
\end{equation}
\end{proposition}
\begin{proof}
Since $\varphi(r)\in[0,\pi/4]$ for all $r\ge e$, we only have to show that for any sequence $(r_n,\varphi_n)$ with $\varphi(r_n)=\varphi_n$ and such that $\varphi_n\to\varphi_0$ and $r_n\to\infty$ we have $4(\log r_n)\sin\varphi_n\to\pi$. 
We know that 
\[4\log r_n\sin\varphi_n\cos\varphi_n=(\pi-2\varphi_n)(1-2\sin^2\varphi_n)+2\sin\varphi_n\cos\varphi_n+\frac{\sin\varphi_n}{r_n}.\]
Here the right-hand side converges to 
\[A=(\pi-2\varphi_0)(1-2\sin^2\varphi_0)+2\sin\varphi_0\cos\varphi_0\in \R.\]
Therefore $4\log r_n\sin\varphi_n\cos\varphi_n\to A$. Since $0\le \varphi_n\le \pi/4$ and 
$\log r_n\to+\infty$, it follows that $\varphi_n\to0$. Hence $A=\pi$.
\end{proof}

\begin{remark}
The function $\varphi(r)$ have an asymptotic expansion for $r\to+\infty$ that starts by 
\begin{equation}\label{limitL}
\varphi(r)=\frac{\pi}{4\log r}-\frac{\pi^3}{48\log^3 r}+\frac{\pi^3}{96\log^4r}+\frac{\pi^5}{320\log^5r}+
\Orden(\log^{-6}r).
\end{equation}
This is obtained making a power series that vanish 
\[(2\log r)\sin(2\varphi)-2(\tfrac{\pi}{2}-\varphi)\cos(2\varphi)-\sin(2\varphi)=0.\]

Also, numerically it is found that 
\begin{equation}
\frac12<\frac{4(\log r)\sin(\varphi(r))}{\pi}\le 1.
\end{equation}
\end{remark}

\begin{theorem}\label{theorem6}
Let $L\subset\C$ defined by 
\begin{equation}\label{D:defL}
L=\{s\in\C\colon |s|>2\pi e^2, -\tfrac{\pi}{2}+2\varphi(\sqrt{|s|/2\pi})<\arg s<\tfrac{\pi}{2}\}.\end{equation}
where $\varphi(r)$ is the function considered in the Lemma \ref{L:region}.

For $s\in L$ with $t>0$ we have
\begin{equation}
\Rzeta(s)=\sum_{n=1}^\ell\frac{1}{n^s}+\Orden(|s/2\pi e|^{-\frac{\sigma}{2}}).
\end{equation}
For $s\in L$ and $t<0$ we have for some constant $c>0$
\begin{equation}\label{secondinequality}
\Rzeta(s)=\sum_{n=1}^\ell\frac{1}{n^s}+\Orden(e^{-c\frac{\sqrt{|s|}}{\log|s|}}).
\end{equation}
\end{theorem}

\begin{proof}
We apply Theorem \ref{Theorem1} with $\theta=\pi/2$ and $K=1$ and we get 
that for $\sigma>0$ and $s\to\infty$
\begin{equation}
\Rzeta(s)=\sum_{n=1}^\ell \frac{1}{n^s}+\frac{(-1)^\ell}{2i}\xi^{-s}e^{\pi i \xi^2}
\Bigl(D_0(q)+\frac{D_1(q)}{\xi}+\Orden(|\xi|^{-2})\Bigr).
\end{equation}
We have $q=q_1+iq_2$ with $-1\le q_1-q_2\le 1$. By Proposition \ref{Gproperties},  $D_0(q)$ and $D_1(q)$  are bounded on this strip.

The hypothesis $\sigma>0$ is equivalent to  $\xi_1>0$ and $\xi_2<0$, i.e.  $\xi$ is in the fourth quadrant. The relationship between $\xi$ and $s$ may be expressed by 
\begin{equation}
\sigma=-4\pi\xi_1\xi_2,\qquad t=2\pi(\xi_1^2-\xi_2^2),
\end{equation}
and a simple computation yields
\begin{equation}
|\xi^{-s}e^{\pi i \xi^2}|=\exp\left(4\pi\xi_1\xi_2\log|\xi|
+2\pi(\xi_1^2-\xi_2^2)\arg(\xi)-2\pi\xi_1\xi_2\right).
\end{equation}

\begin{figure}[H]
\begin{center}
\includegraphics[width=0.6\hsize]{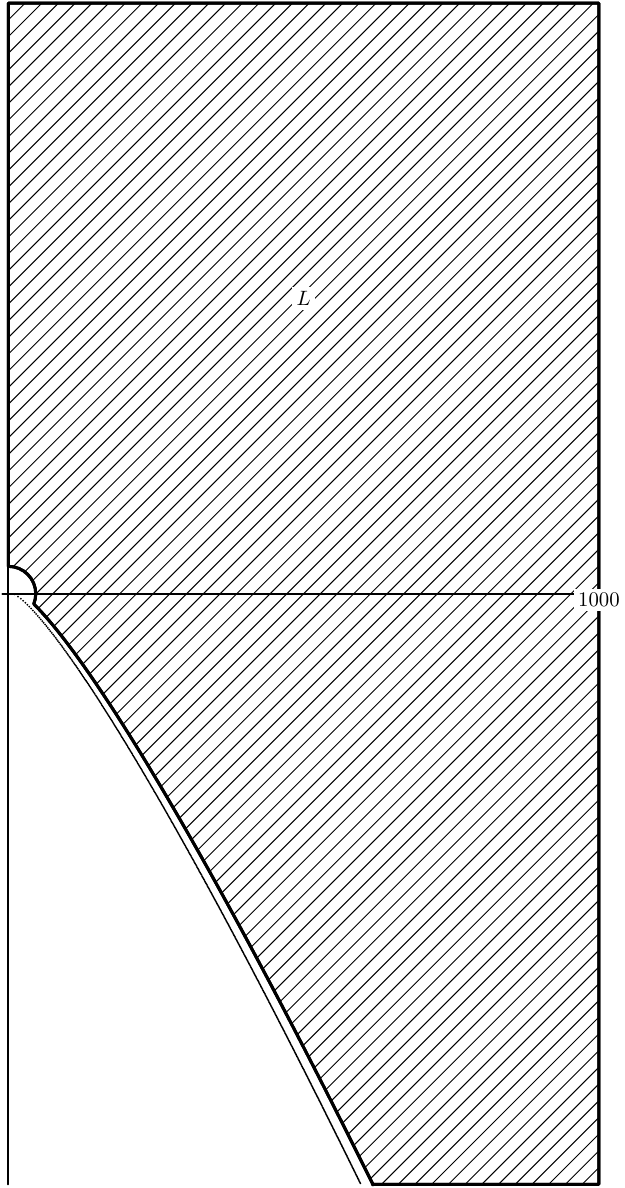}
\caption{The region $L$ in Theorem \ref{theorem6} and the zeros on the fourth quadrant. }
\label{Landzeros}
\end{center}
\end{figure}

When $t>0$, we have  $\xi_1\ge-\xi_2$ the term $2\pi(\xi_1^2-\xi_2^2)\arg(\xi)\le0$ and  for $s\in L$, we have 
\begin{equation}
|\xi^{-s}e^{\pi i \xi^2}|\le e^{4\pi\xi_1\xi_2\log(|\xi|e^{-1/2})} =\Bigl|\frac{s}{2\pi e}\Bigr|^{-\frac{\sigma}{2}}.
\end{equation}
This proves the first assertion.

When  $t<0$ we have $-\pi/2<\arg\xi<-\pi/4$. 
Since $s\in L$ 
\[-\tfrac{\pi}{2}+2\varphi(\sqrt{|s|/2\pi})<\arg s<\tfrac{\pi}{2}.\]
We have $2\arg\xi=\arg s-\frac{\pi}{2}$. Therefore,
\[-\tfrac{\pi}{2}+\varphi(|\xi|)<\arg\xi<0.\]
Since $-\pi/2<\arg\xi<-\pi/4$ and $|\xi|>e$,  we get $\xi\in P$ (see Lemma \ref{L:region}) and it follows that $|\xi^{-2\pi i\xi^2}e^{\pi i\xi^2+\pi\xi}|\le 1$.
Therefore,
\[|\xi^{-2\pi i\xi^2}e^{\pi i\xi^2}|\le |e^{-\pi\xi}|=e^{-\pi\xi_1(s)}.\]
We have $\xi_1(s)=|\xi(s)|\sin\phi$, and by Proposition \ref{orderphi}, we have $\sin\varphi>c/\log |\xi|$ for some constant $c>0$. Noticing that  
\[\frac{|\xi|}{\log|\xi|}\gg \frac{\sqrt{|s|}}{\log|s|}.\]
The inequality in \eqref{secondinequality} is proved. 
\end{proof}

\begin{remark} Siegel restricts the domain to the first quadrant. Figure \ref{Landzeros} represents the region $L$ defined in Theorem \ref{theorem6}. The line to the left of the region in the fourth quadrant is obtained by putting a small dot on each of the 472 zeros of $\Rzeta(s)$ contained in the square $(0,1000)\times(0,-1000)$. 
\end{remark}

\begin{corollary}
Let $L$ be the region defined in Theorem \ref{theorem6}. 
There is a certain $r_0>0$ such that  $s\in L$, 
$\sigma\ge2$ and $|s|\ge r_0$, we have
\begin{equation}
|\Rzeta(s)-1|\le \frac34.
\end{equation}
In particular, $\Rzeta(s)$ do not vanish on this closed set.
\end{corollary}

\begin{proof}
When  $s\in L$ we have $\sigma>0$ and $\xi$ is in the fourth quadrant. Since $|\xi|^2=|s/2\pi|$, the conditions $|s|\to\infty$ and $|\xi|\to\infty$ are equivalent for $s\in L$. Since $\xi$ is in the fourth quadrant, we have 
$\xi_1-\xi_2=|\xi_1|+|\xi_2|\sim |\xi|$.
Therefore, $s\to\infty$ in $L$ implies $\ell\to\infty$ (note the definitions in \eqref{defeta}, \ref{defm}). 

By Theorem \ref{theorem6}
we have for $s\in L$ 
\[\Bigl|\Rzeta(s)-\sum_{n=1}^\ell\frac{1}{n^s}\Bigr|\le\begin{cases}
C|s/2\pi e|^{-\sigma/2} &\text{if $t>0$},\\
Ce^{- c\sqrt{|s|}/\log|s|} &\text{if $t<0$.}\end{cases}\]
Take a $\delta>0$ with  $\pi^2/6-1+\delta<3/4$. Then for $s\in L$, $\sigma>2$ and $|s|$ big enough, we will have
\[\Bigl|\sum_{n=1}^\ell\frac{1}{n^s}-1\Bigr|+\delta\le \delta+\sum_{n=2}^\ell\frac{1}{n^2}
\le \delta +\frac{\pi^2}{6}-1<\frac34.\]
When $t\ge0$ the error $C|s/2\pi e|^{-\sigma/2}\le C|s/2\pi e|^{-1}<\delta$ if $s$ is large enough.  Analogously, when $t<0$ we have $Ce^{-c\sqrt{|s|}/\log|s|}<\delta$ for $s$ large enough.  This proves our claim. 
\end{proof}

\begin{corollary}
For $\sigma\ge1/2$ and $t>0$ we have $\Rzeta(s)=\Orden(t^{1/4})$ when $s\to\infty$.
\end{corollary}

\begin{proof}
We apply Theorem \ref{theorem6}.
For $\sigma\ge 2$ and $t>0$ we have $\sum_{n=1}^\ell\frac{1}{n^s}=\Orden(1)$.  For $\tfrac12 \le\sigma\le 2$ we have $\ell=\Orden(t^{1/2})$ and the sum is $\Orden(t^{1/4})$.
\end{proof}

\section{Asymptotic on the third and fourth quadrants.}

\begin{theorem}\label{asympinf} For a given $0<\theta<\pi/2$ let $M=M_\theta\subset \C$ the closed set
\begin{equation}
M=\{s\in\C\colon|s|>2\pi e^2,\quad -\pi+\theta<\arg s<-\tfrac{\pi}{2}+\arctan
\tfrac{\pi}{2\log|\xi|}\}.
\end{equation}
On this set we have
\begin{equation}\label{mainpart}
\Rzeta(s)=\frac{(-1)^\ell}{2i}\xi^{-s}e^{\pi i \xi^2}G(q)\bigl(1+\Orden(|\xi|^{-1})\bigr).
\end{equation}
\end{theorem}

\begin{proof} 
First, notice that for $s\in M$ we have $|\xi|>e$. So,  for $s\in M$ we have 
$-\pi+\theta<\arg s<0$, and $-\frac{3\pi}{4}+\frac{\theta}{2}<\arg\xi<-\frac{\pi}{2}+\frac12\arctan\frac{\pi}{2\log|\xi|}$.  Also  $|\xi|^2=|s|/2\pi$  and  $s\to\infty$ is equivalent to $\xi\to\infty$. 

Applying Theorem \ref{Theorem1} with $K=1$, we get for all $s\in M$
\begin{equation}
\Rzeta(s)=\sum_{n=1}^\ell n^{-s}+\frac{(-1)^\ell}{2i}\xi^{-s}e^{\pi i \xi^2}\Bigl(D_0(q)+\frac{D_1(q)}{\xi}+\Orden_\theta(e^{-\frac{\pi}{2\sqrt{2}}|\mu|}|\xi|^{-2})\Bigr)
\end{equation}
The explicit formulas for the coefficients $D_k(q)$ in \eqref{explicit}, the inequalities in \eqref{E:Gineq1} and \eqref{E:Gineq2} and the fact that $G(q)$ does not vanish for $q\in B_1$ implies that 
\begin{equation}
\Rzeta(s)=\frac{(-1)^\ell}{2i}\xi^{-s}e^{\pi i \xi^2}G(q)\Bigl(1+\frac{2i(-1)^\ell
\xi^{s}e^{-\pi i \xi^2}}{G(q)}\sum_{n=1}^\ell n^{-s}+\Orden_\theta(|\xi|^{-1})\Bigr).
\end{equation}
We only need to show that for $s\in M$
\[\Bigl|\frac{\xi^{s}e^{-\pi i \xi^2}}{G(q)}\sum_{n=1}^\ell n^{-s}\Bigr|\le \frac{c}{|\xi|}.\]
Let $\varphi=-\frac{\pi}{2}-\arg\xi\in[-\frac12\arctan\frac{\pi}{2\log|\xi|}, \frac{\pi}{4}-\frac{\theta}{2}]$, and $r=|\xi|$.  In terms of $r$ and $\varphi$ with $s=\sigma+it$ we have 
\[|\xi|=r,\quad \xi_1=-r\sin\varphi,\quad \xi_2=-r\cos\varphi,\quad \arg(\xi)=-(\tfrac{\pi}{2}+\varphi).\]
\[\ell=\lfloor r(\cos\varphi-\sin\varphi)\rfloor,\quad \sigma=-4\pi r^2\sin\varphi\cos\varphi,\quad t=-2\pi r^2(\cos^2\varphi-\sin^2\varphi).\]
Notice that in the range of $\varphi$ considered we have $\cos\varphi-\sin\varphi>0$, so that for $r=|\xi|\to+\infty$ we have $\ell\to+\infty$. 
\[q_1=2\xi_1-2\ell-1=-2r\sin\varphi-2\ell-1,\quad q_2=2\xi_2=-2r\cos\varphi,
\]
\begin{align*}
\mu&=\frac{1}{\sqrt{2}}\bigl(-2r(\sin\varphi+\cos\varphi)-2\lfloor r(\cos\varphi-\sin\varphi)\rfloor-1\bigr)\\
&=\frac{1}{\sqrt{2}}\bigl(-4r\cos\varphi+2r(\cos\varphi-\sin\varphi)-2\lfloor r(\cos\varphi-\sin\varphi)\rfloor-1\bigr)\\
&=-2\sqrt{2} r\cos\varphi+w, \quad |w|\le 2^{-1/2}.
\end{align*}

We consider apart the cases $\sigma\le0$  and $\sigma>0$. 

Consider first the case $\sigma\le0$, or equivalently $\varphi\ge0$. In this case
\[\Bigl|\sum_{n=1}^\ell n^{-s}\Bigr|\le \ell \ell^{-\sigma}\le r(\cos\varphi-\sin\varphi)\exp\bigl(2\pi r^2(\log(r(\cos\varphi-\sin\varphi))\sin2\varphi
\bigr)\]
\begin{align*}
|\xi^s e^{-\pi i \xi^2}|&=\exp(-4\pi\xi_1\xi_2\log|\xi|-2\pi(\xi_1^2-\xi_2^2)\arg(\xi)+2\pi\xi_1\xi_2)\\
&=\exp\bigl(-2\pi r^2(\log r)\sin2\varphi-2\pi(\tfrac{\pi}{2}+\varphi) r^2\cos2\varphi+\pi r^2\sin2\varphi\bigr).
\end{align*}
By \eqref{E:Gineq1} 
\[G(q)^{-1}\ll e^{\frac{\pi}{2\sqrt{2}}|\mu|}\ll e^{\pi r\cos\varphi}\le e^{\pi r}.\]
It follows that 
\[\Bigl|\frac{\xi^{s}e^{-\pi i \xi^2}}{G(q)}\sum_{n=1}^\ell n^{-s}\Bigr|\ll r e^{\pi r+\pi r^2 W},\]
where 
\[W:=-2\sin(2\varphi)\log\frac{1}{\cos\varphi-\sin\varphi}-2(\tfrac{\pi}{2}+\varphi) \cos2\varphi+\sin2\varphi\]
This is a continuous function on $[0,\pi/4]$ that is always negative. 
In fact, we have
\begin{align*}
2(\tfrac{\pi}{2}+\varphi) \cos2\varphi>1,\quad 0\le \varphi\le \pi/6,\\
2\sin(2\varphi)\log\frac{1}{\cos\varphi-\sin\varphi}>1,\quad \pi/6<\varphi\le \pi/4.
\end{align*}
Therefore, there exists $W_0>0$ such that 
\[\Bigl|\frac{\xi^{s}e^{-\pi i \xi^2}}{G(q)}\sum_{n=1}^\ell n^{-s}\Bigr|\ll r e^{\pi r-\pi r^2 W_0}.\]

For $\sigma>0$ we change the bound of the zeta sum
\[\Bigl|\sum_{n=1}^\ell n^{-s}\Bigr|\le \ell \le r(\cos\varphi-\sin\varphi).\]
Except that now $\varphi<0$ the above bounds are valid and we obtain
\[\Bigl|\frac{\xi^{s}e^{-\pi i \xi^2}}{G(q)}\sum_{n=1}^\ell n^{-s}\Bigr|\ll r e^{\pi r+\pi r^2 W},\]
with a new $W$ that is given by 
\[W:=-2(\log r)\sin2\varphi-2(\tfrac{\pi}{2}+\varphi) \cos2\varphi+\sin2\varphi\]
with $-\frac12\arctan\frac{\pi}{2\log|\xi|}<\varphi<0$,
so that 
$0<-\tan(2\varphi)<\frac{\pi}{2\log|\xi|}$. Therefore,
\begin{align*}
W&\le (2\log r-1)\frac{\pi}{2\log r}\cos(2\varphi)-(\pi+2\varphi)\cos(2\varphi),\\
&=-\frac{\pi}{2\log r}\cos(2\varphi)-2\varphi\cos(2\varphi)
\end{align*}
Since $0<-\tan2\varphi<\frac{\pi}{2\log r}$ 
\[0<-2\varphi <\arctan \frac{\pi}{2\log r}<\frac{\pi}{2\log r}-\frac13\Bigl(\frac{\pi}{2\log r}\Bigr)^3+\frac15\Bigl(\frac{\pi}{2\log r}\Bigr)^5.\]
\begin{align*}
W&\le -\frac{\pi}{2\log r}\cos(2\varphi)-2\varphi\cos(2\varphi)\\
&\le -\frac13\Bigl(\frac{\pi}{2\log r}\Bigr)^3\cos(2\varphi)+\frac15\Bigl(\frac{\pi}{2\log r}\Bigr)^5\cos(2\varphi)\\
&\le -\frac14\Bigl(\frac{\pi}{2\log r}\Bigr)^3,\qquad r\ge r_0.
\end{align*}
Therefore,
\[\Bigl|\frac{\xi^{s}e^{-\pi i \xi^2}}{G(q)}\sum_{n=1}^\ell n^{-s}\Bigr|\ll r e^{\pi r-\frac{\pi r^2}{4}(\frac{\pi}{2\log r})^3 }\le \frac{c}{r}, \qquad r\ge r_0.\qedhere\]
\end{proof}

\begin{remark}
According to \eqref{limitL} and \eqref{D:defL} a point $s$ is in $L$ when 
\[-\frac{\pi}{2}+\frac{\pi}{2\log r}-\frac{\pi^3}{24\log^3r}+\frac{\pi^3}{48\log^4r}+\cdots<\arg s<\frac{\pi}{2},\]
%The expansion of $\varphi$ here is obtained term by term computing the asymptotic expansion of  $v(x,y)=2x \sin(2y)-2(\tfrac{\pi}{2}-y)\cos(2y)-\sin(2y),$
%at the point \[\frac{\pi}{4x}+\sum_{j=2}^{N-1} {a_j}x^{-j}+A_N x^{-N}\] and computing $A_N$ so that the first term of the expansion vanish. 

while it is in $M$ when 
\[-\pi+\theta<\arg s<-\frac{\pi}{2}+\frac{\pi}{2\log r}-\frac{\pi^3}{24\log^3r}+\frac{\pi^5}{160\log^5r}+\cdots\]
We will see in \cite{A108} that there is a line of zeros in the strip 
\[\frac{\pi}{2\log r}-\frac{\pi^3}{24\log^3r}+\frac{\pi^5}{160\log^5r}+\cdots<\frac{\pi}{2}+\arg s<\frac{\pi}{2\log r}-\frac{\pi^3}{24\log^3r}+\frac{\pi^3}{48\log^4r}+\cdots\]
 
\end{remark}

\begin{corollary}\label{asymptinf}
For $t\to+\infty$ we have
\begin{multline}
\Rzeta(\tfrac12 -it)=-\frac{1}{\sqrt{2}}\Bigl(\frac{t}{2\pi}\Bigr)^{-\frac14}
\exp\Bigl\{\frac{\pi t}{2}-\Bigl(\frac{\pi t}{2}\Bigr)^{\frac12}\Bigr\}\cdot\\
\cdot
\exp\Bigl\{i\Bigl(\frac{t}{2}\log\frac{t}{2\pi}-\frac{t}{2}+\frac{3\pi}{8}\Bigr)\Bigr\}
\cdot\bigl(1+\Orden(t^{-\frac12})\bigr)
\end{multline}
\end{corollary}

\begin{proof}
Applying Theorem \eqref{asympinf} for $s=\frac12-it$ with $t\gg 0$, we obtain
\[\Rzeta(\tfrac12-it)=\frac{(-1)^\ell}{2i}\xi^{-s}e^{\pi i \xi^2}G(q)(1+\Orden(|\xi|^{-1}).\]
We have
\begin{displaymath}
\xi=\sqrt{-\frac{t}{2\pi}-\frac{i}{4\pi}}=-i\sqrt{\frac{t}{2\pi}}\Bigl(1+\frac{i}{2t}
\Bigr)^{1/2}
\end{displaymath}
It follows that the real and imaginary parts of $\xi$ are 
\begin{equation}
\xi_1=\frac{1}{4\sqrt{2\pi t}}+\Orden(t^{-5/2}),\quad \xi_2=-\sqrt{\frac{t}{2\pi}}
+\Orden(t^{-\frac32}).
\end{equation}
and 
\begin{equation}
|\xi|=\frac{1}{\sqrt{2\pi}}(\tfrac14 +t^2)^{\frac14}=\sqrt{\frac{t}{2\pi}}
(1+\Orden(t^{-2})).
\end{equation}
Then 
\begin{equation}
\log|\xi|=\frac12\log\frac{t}{2\pi}+\Orden(t^{-2})
\end{equation}
and
\begin{equation}
\arg(\xi)=-\arctan\frac{-\xi_2}{\xi_1}=-\frac{\pi}{2}+\arctan\frac{\xi_1}{-\xi_2}=-\frac{\pi}{2}+\frac{1}{4t}+\Orden(t^{-3}).
\end{equation}
After some computations
\begin{equation}
-s\log\xi+\pi i\xi^2
=\frac{\pi t}{2}-\frac14\log\frac{t}{2\pi}+\Orden(t^{-2})+
i\Bigl(\frac{t}{2}\log\frac{t}{2\pi}-\frac{t}{2}+\frac{\pi}{4}+\Orden(t^{-1})\Bigr)
\end{equation}
By definition, $\ell=\lfloor\xi_1-\xi_2\rfloor$, it follows that
\begin{equation}
\ell=\Bigl\lfloor \sqrt{\frac{t}{2\pi}}+\frac{1}{4\sqrt{2\pi t}}+\Orden(t^{-3/2})
\Bigr\rfloor.
\end{equation}
Then
\begin{multline*}
q=2\xi-2\ell-1=\frac{2}{4\sqrt{2\pi t}}-2\ell-1+\Orden(t^{-5/2})-2i\sqrt{\frac{t}{2\pi}}+
i\Orden(t^{-3/2}))=\\
=-2\sqrt{2}\frac{1+i}{\sqrt{2}}\sqrt{\frac{t}{2\pi}}+
2\Bigl(\sqrt{\frac{t}{2\pi}}+\frac{1}{4\sqrt{2\pi t}}-\ell\Bigr)-1+\Orden(t^{-3/2})
\end{multline*}
with  $-1\le \lambda<1$ we have then
\begin{equation}
q=-2e^{\frac{\pi i}{4}}\sqrt{\frac{t}{\pi}}+\lambda+\Orden(t^{-3/2})=
-2e^{\frac{\pi i}{4}}\sqrt{\frac{t}{\pi}}\cdot(1+\Orden(t^{-1/2})).
\end{equation}
Then 
\begin{displaymath}
\frac{\pi i}{2}q^2=\frac{\pi i}{2}4i\frac{t}{\pi}(1+\Orden(t^{-1}))=
-2t(1+\Orden(t^{-1}))
\end{displaymath}
and the term 
\begin{equation}
e^{\frac{\pi i}{2}q^2}=e^{-2t}\Orden(1)
\end{equation}
is very small for $t\to+\infty$. 
Continue with the computation of $G(q)$
\begin{multline*}
\cos\pi q=\cos\pi(2\xi-2\ell-1)=-\cos2\pi\xi=-\frac12e^{2\pi i\xi}
(1+e^{-4\pi i \xi})=\\
=
-\frac12 e^{\sqrt{2\pi t}+\Orden(t^{-3/2})}e^{\frac{2\pi i}{4\sqrt{2\pi t}}+
\Orden(t^{-5/2})}(1+\Orden(e^{-2\sqrt{2\pi t}}))
\end{multline*}
so that
\begin{equation}
\cos\pi q=-\frac{1}{2}e^{\sqrt{2\pi t}}\cdot e^{\frac{2\pi i}{4\sqrt{2\pi t}}}
\cdot
(1+\Orden(t^{-3/2})).
\end{equation}
Analogously,
\begin{multline*}
\cos\frac{\pi q}{2}=\cos\Bigl(\pi\xi-\pi \ell-\frac{\pi}{2}\Bigr)=(-1)^\ell\sin\pi\xi=
\frac{(-1)^\ell}{2i}e^{\pi i\xi}(1-e^{-2\pi i\xi})=\\
=\frac{(-1)^\ell}{2i}e^{\frac{\sqrt{2\pi t}}{2}+\Orden(t^{-\frac32})}\cdot 
e^{\frac{\pi i}{4\sqrt{2\pi t}}+\Orden(t^{-\frac52})}(1+\Orden(e^{-\sqrt{2\pi t}}))
\end{multline*}
so that
\begin{equation}
\cos\frac{\pi q}{2}=\frac{(-1)^\ell}{2i}e^{\frac12\sqrt{2\pi t}}\cdot
e^{\frac{\pi i}{4\sqrt{2\pi t}}}\cdot(1+\Orden(t^{-\frac32})).
\end{equation}
It follows that
\begin{multline*}
\frac{\sqrt{2}e^{\pi i/8}\cos\frac{\pi q}{2}-e^{\frac{\pi i}{2}q^2}}
{\cos\pi q}=\frac{(-1)^\ell}{2i}\frac{\sqrt{2}e^{\pi i/8}e^{\frac12\sqrt{2\pi t}}
e^{\frac{\pi i}{4\sqrt{2\pi t}}}\cdot(1+\Orden(t^{-\frac32}))}
{-\frac{1}{2}e^{\sqrt{2\pi t}}\cdot e^{\frac{2\pi i}{4\sqrt{2\pi t}}}
\cdot
(1+\Orden(t^{-3/2}))}=\\
=-\frac{(-1)^\ell}{i}\sqrt{2}e^{\pi i/8}e^{-\frac12\sqrt{2\pi t}}e^{-\frac{\pi i}{4\sqrt{2\pi t}}}
(1+\Orden(t^{-3/2}))
\end{multline*}
Substituting these values into \eqref{mainpart} the two terms $(-1)^\ell$ cancel. The discontinuity in \eqref{mainpart} is mainly apparent. We obtain
\begin{multline*}
\Rzeta(\tfrac12 -it)=-\frac12\Bigl(\frac{t}{2\pi}\Bigr)^{-\frac14}
e^{\frac{\pi t}{2}}\cdot\\ 
\cdot \sqrt{2}e^{\frac{\pi i}{8}} e^{i\left(\frac{t}{2}\log\frac{t}{2\pi}-\frac{t}{2}+\frac{\pi}{4}
\right)}(1+\Orden(t^{-1}))\cdot e^{-\frac12\sqrt{2\pi t}}
e^{-\frac{\pi i}{4\sqrt{2\pi t}}}
(1+\Orden(t^{-1/2})).
\end{multline*}
That is,
\begin{multline}
\Rzeta(\tfrac12 -it)=
-\frac{1}{\sqrt{2}}\Bigl(\frac{t}{2\pi}\Bigr)^{-\frac14}
\exp\Bigl\{\frac{\pi t}{2}-\Bigl(\frac{\pi t}{2}\Bigr)^{\frac12}\Bigr\}
\cdot\\
\cdot\exp\Bigl\{i\Bigl(\frac{t}{2}\log\frac{t}{2\pi}-\frac{t}{2}+\frac{3\pi}{8}\Bigr)\Bigr\}
(1+\Orden(t^{-1/2})).
\end{multline}
\end{proof}
\begin{remark}
The Riemann-Siegel function $Z(t)$ is given by 
\[Z(t)=2\Re\{e^{i\vartheta(t)}\Rzeta(\tfrac12+it)\}\] (see \cite{A166} for a proof).
Taking the complex conjugate,  we have for $t>0$,
\[Z(t)=2\Re\{e^{-i\vartheta(t)}\Rzeta(\tfrac12-it)\},\] since 
\[\vartheta(t)=\frac{t}{2}\log\frac{t}{2\pi}-\frac{t}{2}-\frac{\pi}{8}+\Orden(t^{-1}).\]
It follows that the product
\[e^{-i\vartheta(t)}\Rzeta(\tfrac12-it)=
-\frac{1}{\sqrt{2}}\Bigl(\frac{t}{2\pi}\Bigr)^{-\frac14}
\exp\Bigl\{\frac{\pi t}{2}-\Bigl(\frac{\pi t}{2}\Bigr)^{\frac12}\Bigr\}\cdot i\cdot(1+\Orden(t^{-1/2}))\]
is almost purely imaginary and very large. We know that the real part is $Z(t)$, which is relatively small. Therefore, for $t<0$ our expression for $Z(t)$ is almost useless.  As Siegel knew, this is not the same for $t>0$. 
\end{remark}

\section{Second expansion}
We have \cite{A166}*{eq.~(2)} 
\begin{equation}\label{uno}
\Rzeta(s)=\chi(s)\bigl\{\zeta(1-s)-\overline{\Rzeta}(1-s)\bigr\}
\end{equation}
where $\overline{\Rzeta}(s):=\overline{\Rzeta(\overline{s})}$, replacing it with its asymptotic expansion we get other expressions for $\Rzeta(s)$. 

We introduce some notation similar to \eqref{defeta}.
Fix $0<\theta<\pi$ and consider $1-s\in\Delta$ (where $\Delta$ is the region defined in \eqref{E:Delta}), we define
\begin{equation}\label{defeta-n}
\eta:=\sqrt{\frac{s-1}{2\pi i}},\quad -\frac{\pi}{4}<\arg\eta<\frac{3\pi}{4},
\quad\eta_1:=\Re(\eta),\quad \eta_2:=\Im(\eta),
\end{equation}
in this case $\eta_1+\eta_2\ge0$ and define  
\begin{equation}\label{defm-n}
m:=\lfloor\eta_1+\eta_2\rfloor,\qquad p:=-2(m+\tfrac12 -\eta).
\end{equation}

\begin{theorem}\label{T:asympcomp}
Let $0<\theta<\pi$ be given. Let $\Delta=\Delta_\theta$ the closed set defined in Theorem \ref{Theorem1}.
Let $K\ge1$ be an integer. Then for $1-s\in\Delta$ we have with the notation in \eqref{defeta-n}, \eqref{defm-n}
\begin{equation}\label{firstexpansion2}
\overline{\Rzeta}(1-s)=
\sum_{n=1}^mn^{s-1}-\frac{(-1)^m}{2i}\eta^{s-1} 
e^{-\pi i \eta^2}\Bigl\{\sum_{k=0}^K\frac{\overline{D}_k(p)}{\eta^k}+
\Orden_{K,\theta}(e^{-\frac{\pi}{2\sqrt{2}}|\mu|}|\eta|^{-K-1})\Bigr\}.
\end{equation}
\end{theorem}
\begin{proof}
Theorem \ref{Theorem1} applied to $z=1-\overline s\in\Delta$ yields
\[\Rzeta(z)=\sum_{n=1}^\ell n^{-z}+\frac{(-1)^\ell}{2i}\xi^{-z} 
e^{\pi i \xi^2}\Bigl\{\sum_{k=0}^K\frac{D_k(q)}{\xi^k}+
\Orden_{K,\theta}(e^{-\frac{\pi}{2\sqrt{2}}|\mu|}|\xi|^{-K-1})\Bigr\}.
\]
Taking the complex conjugate $\overline{\Rzeta(z)}=\overline{\Rzeta}(1-s)$. $\xi=\sqrt{\frac{z}{2\pi i}}$ with $-\frac{3\pi}{4}<\arg\xi<\frac{\pi}{4}$ and it follows that $\overline{\xi}=\eta$.
Therefore $\ell=\lfloor \xi_1-\xi_2\rfloor=\lfloor \eta_1+\eta_2\rfloor=m$. 
Then $\overline{q}=\overline{-2(\ell+\frac12-\xi)}=-2(m+\frac12-\eta)=p$. And with
$q=(\mu+i\nu)e^{\pi i/4}\in B_1$ we obtain $p=(\mu-i\nu)e^{-\pi i/4}\in\overline{B}_1$. 
Therefore, we obtain \eqref{firstexpansion2} by conjugating the above equation.
\end{proof}

\begin{theorem}\label{P:intermediate}
Let $0<\theta<\pi$ be given. For $s\in \C$ with  $\theta<\arg(s-1)<2\pi-\theta$ 
we have with the notations in \eqref{defeta-n}, \eqref{defm-n}
\begin{equation}\label{E:interm}
\Rzeta(s)=\chi(s)\Bigl(\zeta(1-s)-\sum_{n=1}^m \frac{1}{n^{1-s}}+\frac{(-1)^m}{2i}
\eta^{s-1}e^{-\pi i \eta^2}\overline{G}(p)(1+\Orden_\theta(|\eta|^{-1}))\Bigr).
\end{equation}
\end{theorem}
\begin{proof}
In \eqref{uno} we apply Theorem \ref{T:asympcomp} with $K=1$. In the resulting equation, the terms appear
\[\overline{D}_0(p)+\frac{\overline{D}_1(p)}{\eta}+
\Orden(e^{-\frac{\pi}{2\sqrt{2}}|\mu|}|\xi|^{-2}).\]
$D_0(z)=G(z)$ and $\overline{p}\in B_1$, so that $\overline{D}_0(p)\ne0$ and we may write
these three terms as 
\[\overline{G}(p)\Bigl(1+\frac{\overline{D}_1(p)}{\overline{G}(p)\eta}+
\overline{G}(p)^{-1}\Orden(e^{-\frac{\pi}{2\sqrt{2}}|\mu|}|\xi|^{-2})\Bigr).\]
The bounds in \eqref{E:Gineq1} and \eqref{E:Gineq2} transform this into 
\[\overline{G}(p)\bigl(1+\Orden_\theta(|\eta|^{-1})\bigr).\]
Substituting on the expression for $\Rzeta(s)$ yields our result. 
\end{proof}

\begin{theorem}\label{left}
Let $G$ be the set of points $s\in\C$ with $|s-1|>2\pi e$ and 
$\frac{\pi}{2}\le\arg(s-1)\le \frac{\pi}{2}+2\arctan\sqrt{\log\frac{|s-1|}{2\pi}}$. 
With the notations of \eqref{defeta} and \eqref{defm} we have 
\begin{equation}\label{E:85-1}
\Rzeta(s)=\chi(s)\Bigl(\zeta(1-s)-\sum_{n=1}^mn^{s-1}+\Orden(1)\Bigr),\qquad
s\in G,\enspace s\to\infty.
\end{equation}
\end{theorem}

\begin{proof}
We apply Theorem \ref{P:intermediate}. Observe also that by Proposition \ref{Gproperties} the function $\overline G(p)$ is bounded by an absolute constant. 
Therefore, we only need to show that $\eta^{s-1}e^{-\pi i \eta^2}$ is bounded in our region $G$. 

In $G$ we have  $0<\arg(\eta)<\pi/2$ so that $\eta_1$ and $\eta_2\ge0$
\begin{equation}\label{etamod}
|\eta^{s-1}e^{-\pi i\eta^2}|=|\eta^{2\pi i \eta^2}e^{-\pi i\eta^2}|=
e^{-2\pi(\eta_1^2-\eta_2^2)\arg(\eta)-4\pi\eta_1\eta_2\log|\eta|+2\pi\eta_1\eta_2}.
\end{equation}
When $\eta_1\ge\eta_2$ and $\log|\eta|>\tfrac12 $, it is clear that 
$|\eta^{s-1}e^{-\pi i\eta^2}|\le 1$.   When $\eta_1<\eta_2$ we 
have $\arg(\eta)<\frac{\eta_2}{\eta_1}$, so that 
\begin{multline*}
-2\pi(\eta_1^2-\eta_2^2)\arg(\eta)-4\pi\eta_1\eta_2\log|\eta|+2\pi\eta_1\eta_2\le \\
\\\le 2\pi(\eta_2^2-\eta_1^2)\frac{\eta_2}{\eta_1}-4\pi\eta_1\eta_2\log|\eta|
+2\pi\eta_1\eta_2= 2\pi\frac{\eta_2}{\eta_1}(\eta_2^2-2\eta_1^2\log|\eta|)
\end{multline*}
is negative for $\frac{\eta_2}{\eta_1}<\sqrt{2\log|\eta|}$. This condition is equivalent to \[\arg\frac{s-1}{2\pi i}<2\arctan\sqrt{\log\frac{|s-1|}{2\pi}}.\qedhere\] 
\end{proof}

\medskip

\begin{remark}
With Theorem \ref{T: thm N} we will see that in most of the region where Theorem \ref{left} applies, the term we have bounded is smaller than the sum of the others. Thus this Theorem makes sense only for $t>0$ and $-1<-\sigma<t^\alpha$, for some $\alpha>0$.
\end{remark}

\begin{corollary}\label{oldcor}
For $\sigma:=\Re(s)<0$ in the set $G$ defined in Theorem \ref{left} and for $s\to\infty$
\begin{equation}
\Rzeta(s)=\chi(s)\Orden(\log|s|)
\end{equation}
\end{corollary}

\begin{proof}
We apply Theorem \ref{left}. For  $\sigma\le 0$ and $|t|\ge2$ it is well known
\cite{E}*{p.~184 Corollary 1} that 
$\zeta(1-s)=\Orden(\log|t|)$. 
For $|t|\le2$ and $\sigma\le -2$ we have $\zeta(1-s)=\Orden(1)$. 
Also, the sum $\sum_{n=1}^m n^{s-1}$ is $\Orden(1)$
for $\sigma\le-1$ and when $-1\le \sigma\le 0$ we have 
\begin{displaymath}
\Bigl|\sum_{n=1}^m n^{s-1}\Bigr|\le \sum_{n=1}^m\frac{1}{n}=\Orden(\log(2+ m)). 
\end{displaymath}
But in this case $0\le m\le |\eta_1|+|\eta_2|\le \sqrt{2}|\eta|=\Orden(|s|^{1/2})$. 
The result follows now from Theorem \ref{left}.
\end{proof}

\begin{corollary}\label{newcor}
Let $\alpha>0$. In the region defined by 
$t>0$ and $0<-\sigma<t^\alpha$ we have 
\begin{equation}\label{missing}
\Rzeta(s)=\chi(s)\Orden(\log|t|)
\end{equation}
\end{corollary}
\begin{proof}
In the region defined by   $t>0$ and $0<-\sigma<t^\alpha$   we may apply Corollary \ref{oldcor}.  But in this region $\log|t|\sim\log|s|$ when $s\to\infty$. 

This corollary proves and extends  the third formula in  Siegel \cite{Siegel}*{eq. (85)}.
\end{proof}

\textbf{Note.} Next corollary \ref{cor85-2} extends the domain of  the second equation  
contained in Siegel \cite{Siegel}*{eq.~(85)}.  

\begin{corollary}\label{cor85-2}
For $t>0$ and $\sigma<0$ such that $1<1-\sigma<t^{1/2}$ and for $t\to+\infty$, 
we have
\begin{equation}\label{E: cor85-2}
\Rzeta(s)=\chi(s)\Orden\Bigl[\Bigl(\frac{t}{2\pi}\Bigr)^{\sigma/2} |\sigma|^{-1}\Bigr],\qquad  t\to+\infty.
\end{equation}
\end{corollary}

\begin{proof}
First, notice that in the region $1<1-\sigma<t^{1/2}$, with $t>2$,  we have
\begin{align*}
\eta&=\Bigl(\frac{t}{2\pi}+i\frac{1-\sigma}{2\pi}\Bigr)^{1/2}=
\Bigl(\frac{t}{2\pi}\Bigr)^{1/2}\Bigl(1+i\frac{1-\sigma}{t}\Bigr)^{1/2}\\
&=\Bigl(\frac{t}{2\pi}\Bigr)^{1/2}\Bigl(1+i\frac{1-\sigma}{2t}+\frac{(1-\sigma)^2}{8t^2}-i\frac{(1-\sigma)^3}{16t^3}+\cdots\Bigr).
\end{align*}
Therefore, since $\frac{(1-\sigma)}{t}<t^{-1/2}$ we obtain
\begin{equation}
\begin{aligned}
\eta_1&=\Bigl(\frac{t}{2\pi}\Bigr)^{1/2}\Bigl(1+\Orden(\tfrac{(1-\sigma)^2}{t^2})\Bigr)=\Bigl(\frac{t}{2\pi}\Bigr)^{1/2}(1+\Orden(t^{-1})),\\ 
\eta_2&=\frac{1-\sigma}{2\sqrt{2\pi t}}\Bigl(1+\Orden(\tfrac{(1-\sigma)^2}{t^2})\Bigr)=\frac{1-\sigma}{2\sqrt{2\pi t}}(1+\Orden(t^{-1})),\\
\eta_1+\eta_2&=\Bigl(\frac{t}{2\pi}\Bigr)^{1/2}(1+\Orden(t^{-1/2})),\quad m=\lfloor\eta_1+\eta_2\rfloor.
\end{aligned}
\end{equation}

We apply Theorem \ref{P:intermediate}.  
Since we assume that $1-\sigma>1$,  Proposition  14 in \cite{A92} applies and  yields with $x=\eta_1+\eta_2$
\[\Bigl|\zeta(1-s)-\sum_{n=1}^mn^{s-1}\Bigr|\le \frac{x^{\sigma}}{|\sigma|}\Bigl(1+\frac{|\sigma|}{x}\Bigr)\]
For $t\to+\infty$ we have $|\sigma|\le 1+ t^{1/2}$ and $x= \eta_1+\eta_2\gg t^{1/2}$ we obtain 
\[\Bigl|\zeta(1-s)-\sum_{n=1}^mn^{s-1}\Bigr|\ll\Bigl(\frac{t}{2\pi}\Bigr)^{\sigma/2}
 |\sigma|^{-1}.\]

By Proposition \ref{Gproperties} the function $G(p)$ is bounded by an absolute constant. It remains to be bounded
\begin{align*}
|\eta^{s-1}e^{-\pi i \eta^2}|&=|e^{(\sigma-1+it)(\log|\eta|+i\arg\eta)-\frac12(\sigma-1+it)}|\\
&=\exp\bigl\{(\sigma-1)(\log|\eta|-\tfrac12)-t\arg\eta\bigr\}.
\end{align*}
We have $\log|\eta|=\frac12\log\frac{t}{2\pi}+\Orden(t^{-1/2})$. According to our hypothesis 
$\sigma=\Orden(t^{1/2})$. We also have $0<\arg\eta\sim \frac{\eta_2}{\eta_1}=\Orden(t^{-1/2})$.
It follows that 
\[|\eta^{s-1}e^{-\pi i \eta^2}|=\Bigl(\frac{t}{2\pi}\Bigr)^{\frac{\sigma-1}{2}}
\exp\Bigl(\frac{1-\sigma}{2}-t\arg\eta\Bigr).\]
We only need to show that 
\[\Bigl(\frac{t}{2\pi}\Bigr)^{-\frac{1}{2}}|\sigma|\exp\Bigl(\frac{1-\sigma}{2}-t\arg\eta\Bigr)\le C\exp\Bigl(\frac{1-\sigma}{2}-t\arg\eta\Bigr)\]
is bounded.

Since $\eta$ is in the first quadrant we have (recall that $\frac{1-\sigma}{t}=\Orden(t^{-1/2})$)
\begin{align*}
\arg\eta&=\arctan \frac{\eta_2}{\eta_1}=\arctan\Bigl(\frac{1-\sigma}{2\sqrt{2\pi t}}
\Bigl(\frac{2\pi}{t}\Bigr)^{1/2}(1+\Orden(t^{-1}))\Bigr)\\&=
\frac{1-\sigma}{2t}(1+\Orden(t^{-1}))+\Orden(t^{-3/2})=\frac{1-\sigma}{2t}+\Orden(t^{-3/2}).
\end{align*}
Therefore $\frac{1-\sigma}{2}-t\arg\eta=\Orden(t^{-3/2})$ and its exponential is bounded. \end{proof}

\begin{remark}
In equation (85) of Siegel's paper, he writes two formulas valid on the range 
$0<-\sigma<t^{\frac37}$, $t>0$ 
\[\Rzeta(s)=\pi^{s-\frac12}\frac{\Gamma(\frac{1-s}{2} )}{\Gamma(\frac{s}{2})}\Orden
\Bigl(\Bigl(\frac{t}{2\pi}\Bigr)^{\frac{\sigma}{2}}|\sigma|^{-1}\Bigr),\qquad 
\Rzeta(s)=\pi^{s-\frac12}\frac{\Gamma(\frac{1-s}{2} )}{\Gamma(\frac{s}{2})}\Orden
(\log t).\]
Notice that for $\sigma<-2/\log\frac{t}{2\pi}$ we have
\[\Bigl(\frac{t}{2\pi}\Bigr)^{\frac{\sigma}{2}}|\sigma|^{-1}<\frac{1}{2e}\log\frac{t}{2\pi}\]
and for $-2/\log\frac{t}{2\pi}<\sigma<0$ the inequality is reversed. Therefore, the result of Siegel is equivalent to 
\[\Rzeta(s)=\begin{cases}
\pi^{s-\frac12}\frac{\Gamma(\frac{1-s}{2} )}{\Gamma(\frac{s}{2})}\Orden
\Bigl(\Bigl(\frac{t}{2\pi}\Bigr)^{\frac{\sigma}{2}}|\sigma|^{-1}\Bigr) & \text{for 
$2/\log\frac{t}{2\pi}<-\sigma<t^{3/7}$},\\ \\
\pi^{s-\frac12}\frac{\Gamma(\frac{1-s}{2} )}{\Gamma(\frac{s}{2})}\Orden
(\log t) & \text{for $2/\log\frac{t}{2\pi}<-\sigma<0$}.
\end{cases}\]
An analogous remark can be formulated to our Corollaries  \ref{newcor} and \ref{cor85-2}.
\end{remark}

\begin{remark}  The three equations \eqref{E:85-1}, \eqref{missing} and  \eqref{E: cor85-2} extend and prove Siegel's equation \cite{Siegel}*{eq.~(85)}. 
\end{remark}

\begin{theorem}\label{T: thm N}
Let $N$ be the set of complex numbers $s=\sigma+it$  with $|s-1|>2\pi e^2$, and $\sigma<1$ and either  $t\le2\pi e$ or $t>2\pi e$ and $\sigma\le 1-8\pi(\frac{t}{2\pi}\log\frac{t}{2\pi})^{1/2}$. 
For $s\in N$:  
\begin{equation}\label{E:left}
\Rzeta(s)=\frac{(-1)^m}{2i}
\chi(s)\eta^{s-1}e^{-\pi i \eta^2}\overline{G}(p)(1+\Orden(|\eta|^{-1})).
\end{equation}
\end{theorem}
\begin{proof}
Take $0<\theta<\pi/4$, for example, in  Theorem \ref{P:intermediate}. Then, for all  $s\in N$, we have \eqref{E:interm}. Then to prove our theorem, we only need to show that for $s\in N$ 
\[H:=\Bigl|\zeta(1-s)-\sum_{n=1}^m \frac{1}{n^{1-s}}\Bigr|\cdot\Bigl|\frac{\eta^{1-s}e^{\pi i \eta^2}}{\overline{G}(p)}
\Bigr|\ll |\eta|^{-1}.\]

For $s\in N$, we have $\eta$ in the first quadrant, we introduce polar coordinates  $\eta=re^{i\varphi}$ with $0\le\varphi \le \frac{\pi}{2}$ and $r>e$. 
\begin{gather*}
\eta_1=r\cos\varphi,\quad \eta_2=r\sin\varphi,\quad m=\lfloor r(\cos\varphi+\sin\varphi)\rfloor\le r(\cos\varphi+\sin\varphi),\\
1-s=-2\pi i\eta^2=4\pi\eta_1\eta_2-2\pi(\eta_1^2-\eta_2^2)i,\quad \log\eta=\log r+i\varphi\\
\sigma-1=-4\pi\eta_1\eta_2=-2\pi r^2\sin(2\varphi),\quad t=2\pi r^2\cos(2\varphi).
\end{gather*}
\begin{align*}
p&=-2m-1+2\eta=(2(\eta_1+\eta_2)-2\lfloor\eta_1+\eta_2\rfloor-1)-2\eta_2(1-i)\\
&=\bigl(\lambda e^{\pi i/4}-2\sqrt{2}\eta_2\bigr)e^{-\pi i/4}=(\mu-i\nu)e^{-\pi i/4},\end{align*}
so that 
\[\mu=-2\sqrt{2}\eta_2+\lambda/\sqrt{2},\qquad |\lambda|\le1.\]

 Since $\sigma<1$, $1-\sigma>0$. We apply Proposition 14 in \cite{A92}
with $x=r(\cos\varphi+\sin\varphi)$ so that $m=\lfloor x\rfloor$. 
The hypothesis of this proposition are satisfied since
\[x=r(\cos\varphi+\sin\varphi)<2\pi r^2=|1-s|,\quad |s|\ge|s-1|-1\ge2\pi e^2-1>2.\]
We obtain
\begin{align*}
\Bigl|\zeta(1-s)-\sum_{n=1}^m \frac{1}{n^{1-s}}\Bigr|&\le 
7|1-s|x^{\sigma-1}=14\pi r^2x^{-2\pi r^2\sin(2\varphi)}
\\
&\le Br^2\exp\bigl(-2\pi r^2\sin(2\varphi)\log(r(\cos\varphi+\sin\varphi))\bigr)
\end{align*}

We compute  $|\eta^{1-s}e^{\pi i \eta^2}|$
\begin{align*}
\Re\{(1-s)&\log\eta+\pi i\eta^2\}=\Re\{-2\pi i\eta^2\log\eta+\pi i\eta^2\}
\\
&=2\pi r^2(\log r)\sin(2\varphi)+2\pi r^2\varphi\cos(2\varphi)-\pi r^2\sin(2\varphi).
\end{align*}
By \eqref{E:Gineq1} we have 
\[|\overline{G}(p)|^{-1}\ll e^{\frac{\pi}{2\sqrt{2}}|\mu|}\ll e^{\pi\eta_2}=e^{\pi r\sin\varphi}
\]

With these inequalities and  \eqref{E:Gineq1} we get
\[H\ll r^2 e^V,\]
where 
\begin{align*}
V&=2\pi r^2(\log r)\sin(2\varphi)+2\pi r^2\varphi\cos(2\varphi)-\pi r^2\sin(2\varphi)\\&+\pi r\sin\varphi-2\pi r^2\sin(2\varphi)\log(r(\cos\varphi+\sin\varphi))\\
&=\pi r^2\Bigl(-2\sin(2\varphi)\log(\cos\varphi+\sin\varphi)+2\varphi\cos(2\varphi)-\sin(2\varphi)+\frac{\sin\varphi}{r}\Bigr)
\end{align*}
The lemma \ref{L:function} then proves
\[H\ll r^2\exp\Bigl(-\pi r^2\frac{4\log r}{\pi r^2}\Bigr)\ll r^{-2}.\]
This proves \eqref{E:left} when $\eta=re^{i\varphi}$ with $\varphi_0\le\varphi\le \pi/2$ and $r\ge e$ where $\varphi_0=\frac{\sqrt{2\log r}}{r}$. 

Assume now that $s\in N$ and let $\eta=re^{i\varphi}$  corresponding to this $s$, we must prove $\varphi_0\le\varphi\le \pi/2$ and $r\ge e$. If $t<0$, then $\pi/4<\varphi<\pi/2$  and $r=(|s-1|/2\pi)^{1/2}\ge e$. 

If we assume that $s\in N$ with $t>0$, we have $r>e$. We have $r$ and $\varphi$ such that  $s-1=2\pi i r^2e^{2i\varphi}$. Then $\sigma-1=-2\pi r^2\sin(2\varphi)$, $t=2\pi r^2\cos(2\varphi)$.

If $0<t<2\pi e$, and $|s-1|\ge2\pi e^2$, then $\arg(s-1)=\pi-\arctan\frac{t}{|s-1|}$.
And then  $\varphi=\frac{\pi}{4}-\frac12\arctan \frac{t}{|s-1|}$, and we want to show 
$\varphi_0\le \varphi$. That is, 
\[\frac{\sqrt{2\log r}}{r}+\frac12\arctan \frac{t}{2\pi r^2}\le\pi/4.\]
This is true since for $r\ge e$ and   $t\le 2\pi e$ we have
\[\frac{\pi}{4}>\frac{\sqrt{2}}{e}+\frac12\arctan \frac{1}{e}\ge \frac{\sqrt{2\log r}}{r}+\frac12\arctan \frac{e}{r^2}\ge \frac{\sqrt{2\log r}}{r}+\frac12\arctan \frac{t}{2\pi r^2}.\]

For $t>2\pi e$ and $r\ge e$. In this case  we have $t<2\pi r^2$, and since $s\in N$, we have $1-\sigma>8\pi(\frac{t}{2\pi}\log\frac{t}{2\pi})^{1/2}$. Therefore,
\begin{align*}
\varphi&=\frac12\arg(s-1)-\frac{\pi}{4}=\frac12\arctan\frac{1-\sigma}{t}\ge\frac12\arctan \frac{8\pi(\frac{t}{2\pi}\log\frac{t}{2\pi})^{1/2}}{t}\\&=
\frac12\arctan\frac{4(\log\frac{t}{2\pi})^{1/2}}{(\frac{t}{2\pi})^{1/2}}
\ge \frac12\arctan\frac{4\sqrt{2\log r}}{r}\ge \frac{\sqrt{2\log r}}{r}.\qedhere
\end{align*}

\end{proof}

\begin{remark}
Siegel proves that we can enlarge the region $N$ where \eqref{E:left} is valid taking 
$1-\sigma > t^{3/7}$. He also asserts that he can prove this for $1-\sigma > t^{\varepsilon}$
with $\varepsilon>0$ arbitrary. In \cite{A98} we have shown that Siegel's proof does not extend beyond the exponent $2/5$.  It is an open problem to determine the right-hand boundary of the zeros of $\Rzeta(s)$ in the upper half-plane.
\end{remark}

\section{Riemann-Siegel formula.}

From \eqref{uno} and the functional equation, we get 
\begin{equation}
\zeta(s)=\Rzeta(s)+\chi(s)\overline{\Rzeta}(1-s).
\end{equation}
Now we take $0<\theta<\pi$ and apply Theorem \ref{Theorem1}  and Theorem 
\ref{T:asympcomp}. Let $G$ be the closed set defined by $t\ge2\pi$ and 
$\theta\le\arg(s-1)\le\pi$ and $0\le\arg(s)\le\pi-\theta$.  $G$ is essentially an angle 
of $\pi-2\theta$ radians. For $s\in G$  we may apply the two expansions to obtain
\begin{multline}
\zeta(s)=\sum_{n=1}^\ell \frac{1}{n^s}+\frac{(-1)^\ell}{2i}\xi^{-s}e^{\pi i \xi^2}
\Bigl\{\sum_{k=0}^K \frac{D_k(q)}{\xi^k}+\Orden_{K,\theta}(|\xi|^{-K-1})\Bigr\}+\\
+
\chi(s)\sum_{n=1}^mn^{s-1}-\frac{(-1)^m}{2i}\chi(s)\eta^{s-1} 
e^{-\pi i \eta^2}\Bigl\{\sum_{k=0}^K\frac{\overline{D}_k(p)}{\eta^k}+
\Orden_{K,\theta}(|\eta|^{-K-1})\Bigr\}
\end{multline}

\section*{Appendix: Some Lemmas}

\begin{lemma}\label{L:190615-1}
For $a>0$, there is a constant $C(a)>0$ such that for $\lambda\to\infty$
\begin{equation}
I(a,\lambda):=\int_{-\infty}^\infty e^{-a x^2}e^{-|\lambda-x|}\,dx\le C(a)e^{-|\lambda|}.
\end{equation}
\end{lemma}
\begin{proof}
Since $I(a,\lambda)=I(a,-\lambda)$ we may assume that $\lambda>0$. Then we have for $\lambda>1$
\begin{align*}
I(a,\lambda)&=\int_{-\infty}^\lambda e^{-a x^2}e^{-\lambda+x}\,dx+\int_{\lambda}^\infty
e^{-a x^2}e^{-x+\lambda}\,dx\\
&\le e^{-\lambda}
\int_{-\infty}^\infty e^{-a x^2}e^{x}\,dx+\int_{\lambda}^\infty e^{-a x^2}\,dx\\
&\le \frac{\sqrt{\pi}}{\sqrt{a}}e^{1/4a}e^{-\lambda}+\int_{\lambda}^\infty \frac{2 a x}{2a\lambda}e^{-a x^2}\,dx =\frac{\sqrt{\pi}}{\sqrt{a}}e^{1/4a}e^{-\lambda}+\frac{e^{-a\lambda^2}}{2a\lambda}\\
&\le 
\Bigl(\frac{\sqrt{\pi}}{\sqrt{a}}e^{1/4a}+\frac{e^{-a\lambda^2+\lambda}}{2a}\Bigr)e^{-\lambda}
\le \Bigl(\frac{\sqrt{\pi}}{\sqrt{a}}e^{1/4a}+\frac{e^{1/4a}}{2a}\Bigr)e^{-\lambda}.
\end{align*}
Since $I(a,\lambda)$ is continuous in $\lambda$ it is clear that for $0\le \lambda\le 1$ the inequality is true, taking $C(a)$ sufficiently great. 
\end{proof}

\begin{lemma}\label{L:simplebound}
Given positive real numbers $a$, $b$, and $c$, there exists a constant $C=C(a,b,c)$ such that 
\begin{equation}
J(a,b,c,\lambda):=\int_{-\infty}^\infty e^{-ax^2+b|x|}e^{-c|\lambda-x|}\,dx\le C(a,b,c)e^{-c|\lambda|}.
\end{equation}
\end{lemma}
\begin{proof}
Again $J(a,b,c,\lambda)=J(a,b,c,-\lambda)$ and we may assume that $\lambda>1$.  We have 
$e^{b|x|}\le e^{bx}+e^{-bx}$. Therefore, it suffices to prove that for $\lambda>1$ there exist $C_1(a,b,c)$ and $C_2(a,b,c)$ such that 
\begin{align*}
&\int_{-\infty}^\infty e^{-ax^2+bx}e^{-c|\lambda-x|}\,dx\le C_1(a,b,c)e^{-c|\lambda|},\\
&\int_{-\infty}^\infty e^{-ax^2-bx}e^{-c|\lambda-x|}\,dx\le C_2(a,b,c)e^{-c|\lambda|}.
\end{align*}
Consider, for example, the first one, the other is similar. We have
\begin{align*}
&\int_{-\infty}^\infty e^{-ax^2+bx}e^{-c|\lambda-x|}\,dx=
e^{b^2/4a^2}\int_{-\infty}^\infty e^{-a(x-\frac{b}{2a})^2}e^{-c|\lambda-x|}\,dx\\
&=
e^{b^2/4a^2}\int_{-\infty}^\infty e^{-ay^2}e^{-c|\lambda-\frac{b}{2a}-y|}\,dy
=e^{b^2/4a^2}\int_{-\infty}^\infty e^{-ay^2}e^{-|c(\lambda-\frac{b}{2a})-cy|}\,dy\\&=
e^{b^2/4a^2}\frac{1}{c}\int_{-\infty}^\infty e^{-ay^2/c^2}e^{-|c(\lambda-\frac{b}{2a})-y|}\,dy
\end{align*}
Applying Lemma \ref{L:190615-1} we obtain
\[\le e^{b^2/4a^2}\frac{1}{c} C(a/c^2)e^{-c|\lambda-\frac{b}{2a}|}\le 
e^{b^2/4a^2}\frac{1}{c} C(a/c^2)e^{cb/2a}e^{-c|\lambda|}.\]
For $|\lambda|\le1$ we proceed as in Lemma \ref{L:190615-1}.
\end{proof}

\begin{lemma}\label{L:function}
For $r>e$ the function 
\[f(r,\varphi):=\frac{4\log r}{\pi r^2}-2\sin(2\varphi)\log(\cos\varphi+\sin\varphi)+2\varphi\cos(2\varphi)-\sin(2\varphi)+\frac{\sin\varphi}{r}\]
is negative for $\varphi\in[\frac{\sqrt{2\log r}}{r},\frac{\pi}{2}]$. 
\end{lemma}

\begin{proof}
Let $\varphi_0=\sqrt{2\log r}/r$, 
first we prove that $f(r, \varphi_0)<0$ and then that $\partial f/\partial\varphi<0$ for $\varphi\in[\varphi_0,\pi/4]$. Then we prove 
$f(r,\varphi)<0$ for $\varphi\in[\pi/4,\pi/2]$, always assuming that $r>e$. 

The function  $2\varphi\cos(2\varphi)-\sin(2\varphi)$ is decreasing for $\varphi\in[0,\pi/4]$
and vanishes at $\varphi=0$; and $\sin\varphi_0\le\varphi_0$ so that 
\[f(r,\varphi_0)\le \frac{4\log r}{\pi r^2}+\frac{\sqrt{2\log r}}{r^2}-2\sin(2\varphi_0)\log(\cos\varphi_0+\sin\varphi_0).\]
We have the Taylor expansion 
\[-2\sin(2\varphi)\log(\cos\varphi+\sin\varphi)=-4\varphi^2+4\varphi^3-\frac{64}{45}\varphi^6+\cdots\] We  may show that for $r>e$
\begin{align*}
f(r,\varphi_0)&\le \frac{4\log r}{\pi r^2}+\frac{\sqrt{2\log r}}{r^2}-2\sin(2\varphi_0)\log(\cos\varphi_0+\sin\varphi_0)\\
&\le \frac{4\log r}{\pi r^2}+\frac{\sqrt{2\log r}}{r^2}-\frac{8\log r}{r^2}+\frac{4(2\log r)^{3/2}}{r^3}\\
&\le \Bigl(\frac{4}{\pi}+\frac{\sqrt{2}}{\sqrt{\log r}}-8+\frac{8\sqrt{2\log r}}{r}\Bigr)\frac{\log r}{r^2}<0.
\end{align*}

Consider now the derivative on $[\varphi_0,\pi/4]$
\begin{multline*}
\frac{\partial f(r,\varphi)}{\partial \varphi}=
\frac{\cos\varphi}{r}-4\cos(2\varphi)\log(\cos\varphi+\sin\varphi)-4\varphi\sin(2\varphi)\\-
\frac{2(\cos\varphi-\sin\varphi)\sin(2\varphi)}{\cos\varphi+\sin\varphi}.
\end{multline*}
For $0<\varphi<\frac{\pi}{4}$ all terms are negative except for the first. In the interval $\varphi_0\le \varphi\le \pi/6$
\begin{align*}
\frac{\cos\varphi}{r}-\frac{2(\cos\varphi-\sin\varphi)\sin(2\varphi)}{\cos\varphi+\sin\varphi}&=\frac{\cos\varphi}{r}-\frac{2(\cos\varphi-\sin\varphi)\sin(2\varphi)}{(\cos\varphi+\sin\varphi)\varphi}\varphi\\
&\le \frac{1}{r}-\frac{18-6\sqrt{3}}{\pi(1+\sqrt{3})}\frac{\sqrt{2\log r}}{r}<0, \qquad r>e.
\end{align*}
For $\frac{\pi}{6}\le \varphi\le \frac{\pi}{4}$ we have
\[\frac{\cos\varphi}{r}-4\varphi\sin(2\varphi)\le \frac{1}{r}-\frac{\pi}{\sqrt{3}}<0.\]

Only remains to show that $f(r,\varphi)<0$ for $\pi/4\le\varphi\le\pi/2$.
In the interval $\frac{\pi}{4}\le \varphi\le \frac{\pi}{2}$ the term $-2\sin(2\varphi)\log(\cos\varphi+\sin\varphi)\le 0$ and therefore
\[f(r,\varphi)\le \frac{4\log r}{\pi r^2}+\frac{\sin\varphi}{r}+2\varphi\cos(2\varphi)-\sin(2\varphi).\]
Since $\cos(2\varphi)<0$ we have 
\[f(r,\varphi)\le \frac{4\log r}{\pi r^2}+\frac{\sin\varphi}{r}+\cos(2\varphi)-\sin(2\varphi)\le \frac{4\log r}{\pi r^2}+\frac{1}{r} -1<0, \qquad r>e.\qedhere\]
\end{proof}

\end{document}